\documentclass[a4paper,11pt,english]{article} 
\usepackage{amsmath}
\usepackage{amsfonts}
\usepackage{amssymb}
\usepackage[all]{xypic}
\usepackage{graphics}
\usepackage{babel}
\usepackage{hyperref}
\usepackage{here}
\usepackage{verbatim}

\def\a{\alpha}
\def\b{\beta}

\def\e{\varepsilon}

\def\vp{\varphi}

\def\o{\omega}

\def\bR{{\mathbb R}}
\def\bZ{{\mathbb Z}}
\def\bC{{\mathbb C}}

\def\b1{{\rm id}}

\newfont{\goth}{eufm10 scaled \magstep1}
\def\ga{\mbox{\goth a}}

\def\gc{\mbox{\goth c}}

\def\gsu{\mbox{\goth su}}

\def\gg{\mbox{\goth g}}
\def\gh{{\mbox{\goth h}}}
\def\gk{\mbox{\goth k}}
\def\gl{\mbox{\goth l}}
\def\gm{\mbox{\goth m}}
\def\gn{\mbox{\goth n}}

\def\gp{\mbox{\goth p}}
\def\gq{\mbox{\goth q}}

\def\gt{\mbox{\goth t}}
\def\gu{\mbox{\goth u}}

\newfont{\mcal}{eusm10 scaled \magstep1}

\def\ch{\mbox{\mcal H}}

\newtheorem{Th}{Theorem}
\newtheorem{theorem}{Theorem}
\newtheorem{proposition}[theorem]{Proposition}
\newtheorem{definition}[theorem]{Definition}

\newtheorem{remark}[theorem]{Remark}
\newtheorem{remar}[theorem]{Remark}
\newtheorem{Prop}[Th]{Proposition}
\newtheorem{Cor}[Th]{Corollary}
\newtheorem{Lem}[Th]{Lemma}
\newtheorem{Def}[Th]{Definition}
\newtheorem{Ex}[Th]{Example }
\def\bt{\begin{Th}}
\def\et{\end{Th}}
\def\bp{\begin{Prop}}
\def\ep{\end{Prop}}
\def\bc{\begin{Cor}}
\def\ec{\end{Cor}}
\def\bl{\begin{Lem}}
\def\el{\end{Lem}}
\def\bd{\begin{Def}}
\def\ed{\end{Def}}
\def\bex{\begin{Ex}}
\def\eex{\end{Ex}}
\def\br{\begin{remar}}
\def\er{\end{remar}}

\def\qed{\hfill$\square$}

\def\be{\begin{equation}}
\def\ee{\end{equation}}
\def\ben{\begin{enumerate}}
 \def\een{\end{enumerate}}
\def\ba{\begin{array}{rlll}}
\def\ea{\end{array}}
\def\bea{\begin{eqnarray}}
\def\eea{\end{eqnarray}}
\def\bean{\begin{eqnarray*}}
\def\eean{\end{eqnarray*}}

\def\ker{\mathrm{ker\;}}

\def\C1{cohomogeneity one}

\begin{document}

\title{Cohomogeneity one K\"ahler  and K\"ahler-Einstein manifolds  with   one  singular orbit II }
\author{Dmitri Alekseevsky\footnote{A.A.Kharkevich Institute for Information Transition
Problems, B.Karetnuj  per.,19, 127051, Moscow, Russia and
Faculty of Science  University of Hradec Kralove,  Rokitanskeho 62, Hradec Kralove,
50003, Czech Republic } \and  Fabio Zuddas\footnote{Dipartimento di Matematica e Informatica, Via Ospedale 72, Cagliari (Italy)}}



\maketitle

 \abstract
  F. Podest\`a  and A. Spiro \cite{P-S}  introduced   a   class  of  $G$-manifolds $M$  with
      a cohomogeneity  one  action   of  a  compact  semisimple  Lie  group $G$   which  admit an  invariant
       K\"ahler   structure  $(g,J)$ (``standard  $G$-manifolds")   and   studied  invariant  K\"ahler  and  K\"ahler-Einstein  metrics  on  $M$.\\
 In  the  first part of this  paper,  we    gave  a  combinatoric   description   of  the   standard  non  compact  $G$-manifolds      as   the  total  space
   $M_{\varphi}$  of  the   homogeneous vector  bundle $M =  G\times_H  V \to  S_0 =G/H$  over  a    flag manifold  $S_0$ and we gave  necessary  and  sufficient  conditions   for the existence of an  invariant K\"ahler-Einstein metric  $g$   on   such manifolds  $M$ in  terms  of the  existence  of  an interval  in  the  $T$-Weyl  chamber  of  the flag manifold $F =  G \times _H PV$  which  satisfies  some  linear  condition.\\
  In  this  paper,   we  consider   standard  cohomogeneity one  manifolds  of  a  classical  simply  connected  Lie  group  $G =  SU_n, Sp_n. Spin_n$   and
    reformulate    these  necessary  and  sufficient  conditions in  terms  of easily   checked    arithmetic properties  of  the  Koszul numbers  associated  with  the  flag manifold $S_0 = G/H$. If  this  conditions  is  fulfilled,  the  explicit  construction  of the  K\"ahler-Einstein metric  reduces  to the calculation  of  the  inverse  function  to a  given function  of  one  variable.

\tableofcontents

\section{Introduction}

F. Podest\`a  and A. Spiro  \cite{P-S}  defined  a   class of  cohomogeneity  one  $G$- manifolds  $M$, called   standard,
  of   a semisimple  compact Lie  group  $G$   with an invariant  complex  structure $J$.  It is  defined  by
the    condition  that  the  complex  structure  $J$  restricted  to a  regular orbit  $G/L$ defines   a projectable
 CR  structure  $(\mathcal{H},J^{\mathcal{H}})$, so that   the restriction  $\mu:
 {G/L}  \to F =  G/K $  of  the  moment  map to $G/L$
 is  a  holomorphic  map    to    a  flag manifold   $F =  G/K =  G/N_G(L)$   with  a  fixed  invariant  complex  structure $J^F$
 which  does not  depend  on  the  regular orbit  $G/L$.  They  gave  a  nice   description  of the invariant K\"ahler metrics  on  the    complex  manifold
 $(M, J)$   in  terms  of  some  interval  in    the  $T$-Weyl  chamber   ossociated  with  the   complex  structure  $J^F$    and   in  the  case  when
 $M$ is  compact  (hence, it has  two  singular orbits) got necessary  and  sufficient  conditions for  $M$ to admit  an invariant  K\"ahler-Einstein metric.
 Similar  results had  been  obtained    by  A. Dancer  and M.Y. Wang \cite{D-W},  who used  a  different  approach.\\

 In  the  previous paper we  showed     that non  compact  standard  cohomogeneity one manifolds  are exactly  the  total  spaces  $M_{\varphi}$  of      the   homogeneous  complex  vector  bundles    $M_{\varphi} = G \times_H  V_{\varphi}  \to  S_0 = G/H$ over    a  flag manifold   $S_0$  with  an invariant  complex  structure $J^S$   defined   by a representation  $\varphi : H \to   GL(V_{\varphi}), \,  V_{\varphi}  =  \mathbb{C}^m$    with   $\varphi(H)  =  U(V_{\varphi}) \simeq  U_m$   and  gave  a    description  of the invariant   K\"ahler  structures  in   terms of  the  painted  Dynkin  diagrams  associated  with  the  flag manifolds
  $S_0 =  G/H$  (the  singular orbit)  and  $F = G/K = \mu(G/L)  =   G \times_H PV_{\varphi}$  (the projectivisation  of  the  vector  bundle  $M_{\varphi}$). We  also   gave  necessary  and  sufficient  conditions  (similar  to  the  conditions  by  Podest\`a-Spiro)    for the  existence  of   invariant  K\"ahler-Einstein metrics in terms  of an interval in   the  $T$-Weyl   chamber  associated  with  the   complex  structure  $J^F$.  If  this  condition  is   satisfied,  the construction  of an  associated K\"ahler-Einstein metric  is  described explicitly  in terms of  a function $f(t)$  which is   the  inverse  function   to a  function $ t = t(f)$  given by the  integral of an  explicit  function of one  variable.\\

   In the present paper, for  a  non  compact  standard cohomogeneity one
 $G$-manifold $M_{\varphi}  =  G \times_H V_{\varphi}$  of  a  classical   semisimple  Lie  group  $SU_n, Sp_n,  Spin_n$, we reformulate   the necessary and sufficient conditions for the existence of invariant K\"ahler-
Einstein  metrics  on  $M_{\varphi}$  in terms of easily   checked    arithmetic properties  of  the  Koszul numbers  associated
with  the  flag manifold $S_0 =  G/H$,   see Theorems \ref{prop triple2} and \ref{prop triple4}. \\

We  will  always  assume   that  the  group  $G$
is simply connected     and it  acts    on   $M$    almost  effectively.\\

\begin{remar}
\rm

When  the  paper  was  finished  we  find the two interesting  papers \cite{AZ-B} and
\cite{C}, where  invariant Ricci-flat metrics  on  some  holomorphic bundles over  flag manifold  are  constructed.\\
In   \cite{C}   the   author    gets   a  nice  general formula
   for   the unique  asymptotically  conical  Ricci-flat  K\"ahler metric on  the canonical  bundle   $K_F$  of  a  flag  manifolds    $F = G^{\mathbb{C}}/P$ .
In    \cite{AZ-B},   the  authors  describe  more explicitly  the  invariant   Ricci-flat K\"ahler 
metric  on  the  canonical bundle  $K_F$ of  the  reducible  flag manifold
$F =  SU_n/S(U_{n_1} \times \cdots \times U_{n_s})  \times SU_{q}/U_{q-1}$  and  show that,  in  the  case  when   the     $q$-root    $K_F^{\frac{1}{q}}$  exists,  the same  formula  gives  a Ricci  flat Kahler  metric  on the  rank  $q$  holomorphic  vector  bundle
$q   K_F^{\frac{1}{q}}$ .
\end{remar}

\section{Preliminary  and statement of the main results}
\subsection{Cohomogeneity  one  K\"ahler manifolds of  standard  type }
\noindent Following  \cite{P-S}, we     focus  our attention to   \C1 K\"ahler $G$-manifolds $(M, J, \o)$  of  the {\it standard  type}, i.e.  manifolds   which satisfy  the  following   conditions:\\
(i)   a regular  orbit $ S= G x = G/L$ is  an {\it ordinary} manifold. This
means that the  normalizer   $K = N_G(L)$ of  the  stability
subgroup $L$ is  the centralizer  of a  torus in $G$ and $\dim K/L =1$.\\
 In  particular, $F= G/K$ is   a  flag manifold   with  induced    invariant  complex  structure  $J^F$. \\
 (ii) the CR   structure induced  by the complex
structure $J$ of $M$  on  a regular orbit  $S = G/L$  is {\it projectable},
that is  the  restriction  $\pi: S=G/L \to F=G/K$ to   $S$  of  the momentum  map   is a holomorphic map
of a CR manifold  onto the flag manifold  $F = G/K$ equipped  with a
fixed invariant complex  structure $J^F$  (which  does not  depend on $S$).\\
(iii) The  $G$- manifold  $M$   has  only  one  singular  orbit  $S_0 = G/H$,  which is  a  complex  submanifold,  hence $M$   is  not  compact.\\
Condition (ii) depends on the complex  structure  $J$ on $M$    and shows   that the  CR   structure on  a  regular orbit  $G/L$ is  determined  by the invariant  complex  structure $J_F$ on   the  flag manifold $F$.  In particular, all  regular orbits  are isomorphic  as   homogeneous   CR manifolds.\\

\smallskip
 Such a \C1 K\"ahler $G$-manifold $(M, J, \o)$ is  called , shortly,  a  {\it standard  \C1 manifold}.\\

In \cite{AZ} we have proved that any   standard   \C1  manifold  $M$ is the  total
space  of the homogeneous  vector  bundle  (called  {\it admissible  bundle})
$$\pi: M_{\vp} = G\times_H
V_{\vp}\to S_0 =G/H$$
 over  the  singular  orbit $S_0$ defined by
 a  representation $\varphi : H \to  GL(V_{\vp})$  with    $\vp(H) \simeq U_m$ ($m = \dim(V_{\vp})$),   called  {\it admissible representation}.\\

Note  that  the   singular orbit $S_0$  is  the  zero  section  of  the vector  bundle  $\pi$ and   all the other orbits have the  form
 $S_t =  G (tv_0)  = G/L$ and are   regular,   with $   0 \neq v_0 \in V$ a  fixed  vector   and  $t >0$. So  the  set of  $G$-regular  points is    $M_{reg} :=   M_{\varphi} \setminus  S_0  =  G/L \times  \bR^+$,  where   $L \subset H$ is  the  stabilizer  of    the  ray    $ \bR^+ v_0$.\\
It is known (see \cite{P-S}, \cite{AZ})   that   the  singular orbit $S_0 =G/H  \subset  M_{\vp}$  is  a  complex  submanifold,  hence  a  flag manifold. The  induced   complex  structure  $J^S$  on $S_0$  together   with  one   of the $\vp(H) = U(V_{\vp})$-invariant  complex  structures  $\pm J^{V_{\vp}}$  on $V_{\vp}$   defines   the invariant  complex  structure  on  the  manifold  $M_{\vp}$.\\

 \subsection{Examples of    standard \C1   K\"ahler  and  K\"ahler-Einstein manifolds }
  Let   $(F = G/K,  g^F, \omega^F, J^F )$ be  a   homogeneous  K\"ahler manifold   of  a semisimple  compact Lie  group  $G$  with integral
  K\"ahler  form  $\omega$   (a Hodge manifold). Denote  by  $\omega  \in \Lambda^2 (\gg)^*$  the    form defined  by  $\omega^F$ in  the  Lie  algebra $\gg$. It is exact, i.e. $\omega  =  d \sigma, \,  \sigma  \in  \gg^*$. We set  $Z = B^{-1} \sigma$. Then  $F$  is identified  with  the  coadjoint  orbit  of $\sigma$  and  the  adjoint orbit of   $Z$:          $ F =\mathrm{Ad}_G^*\sigma = \mathrm{Ad}_G Z = G/K $.
  Denote  by  $\gl$  the  kernel of $\sigma$ in  $\gk$. Then $\gk = \gl + \bR Z$   is  a $B$-orthogonal  decomposition   and $\gl$ generates  a  closed  subgroup $L$ of  $K$  such  that   $ \pi:S = G/L \to F = G/K$  is a principal $T^1$ bundle   and   $S = G/L$   is  an ordinary manifold. The   form $\sigma$  is   $\mathrm{Ad}_K$-invariant   and  it  extends  to  an  invariant contact  form  $\sigma^S$  on  $S$  which is   a   connection  form of  a   $G$- invariant   connection  with  curvature $\omega^F$.   The complex  structure $J^F$  defines a projectable invariant  CR  structure $J^S$ in the contact  distribution  $\ch = \mathrm{ker \sigma^S}$.
  It is  known  (see \cite{A-C-H=K}, Theorem 2.3)  that   the   K\"ahler metric $g^F$ is extended  to  an invariant  Sasaki metric   $g^S =  \sigma^2 + \frac12 \pi^* g^F$   on  $S$.
  The    $G$-invariant  extension $Z^S$ of     the  $\mathrm{Ad}_K$-invariant vector   $Z$  is  the fundamental vector  field    of the principal bundle  $\pi$  and  it is a Killing   vector  field   for  $g^S$.\\
  The Riemannian  cone $(M = C(F) :=  \bR^+  \times S, g =  dr^2  + r^2 g^S) $  is  a \C1  K\"ahler  $G$-manifold with complex  structure defined  by   (see  \cite{A-C-H=K},  Theorem  2.8)
  $$ J|_{\ch} =  J^S, \,   J \xi = Z^S . $$
where  $\xi =  r \partial_r$ is  the   homothetic  vector  field.
Moreover, if $ (F = G/K, g^F, J^F )$ is  a K\"ahler-Einstein homogeneous manifold,  then $ (S= G/L, g^S, Z^S )  $  is  a Sasaki-Einstein  homogeneous  manifold   and
  the  cone   $(M = C(S), g, J)$ is   a Ricci  flat  K\"ahler \C1 manifold  (see  \cite{B-G}, \cite{Sp}).\\

Note  that   the  cone manifold $M$  is  a \C1 $G$-manifold,  but  it  admits  a transitive  group of   homothetic  transformations, generated  by  $G$  and   the  1-parameter homothety  group $\mathrm{exp} \bR \xi$.\\
We  give  a  generalisation of  this  construction of  K\"ahler  cones   associated  to a  homogeneous  K\"ahler manifold.

 \subsection{Description  of  admissible   vector  bundles}
  We   recall the description  of  the  admissible vector  bundles  $M_{\vp } \to S_0 = G/H$  of rank  $m$ over  a  given   flag manifold  $S_0 = G/H$,  see  \cite{AZ}  for  details.\\

 A  flag manifold   $S_0 = G/H$   is  described  by  a painted  Dynkin  diagram,  which represents a decomposition $\Pi = \Pi_B \cup \Pi_W$ of the system  $\Pi$ of  simple  roots of  $G$ into the subsystem  of  {\it white roots}  $\Pi_W$,  which  corresponds to  the  semisimple  part  $\mathfrak{h}'$  of  the   stability  subalgebra  $\mathfrak{h}$, and the subsystem of {\it black roots} $\Pi_B$. Associated  with black  roots $\beta_i$  fundamental  weights  $\pi_i$  define  a  basis   $B^{-1} \pi_i $ of  the  center  $Z(\mathfrak{h})$, where $B$ is the Killing form, see Appendix  for  details.\\
 Now  we  give  a  short  description  of the admissible  vector  bundles $M_{\vp }= G\times_H V_{\vp} \to  S_0 =G/H$  of  rank $m = \dim(V_{\vp})$ over  a  flag manifold .\\

 \subsubsection{Case of  line  bundles }

 If $m=1$, then   $M_{\vp} =  G\times_H \mathbb{C}$ is   a  complex  line  bundle  defined   by  a character $\chi : Z(H) = T^k \to T^1 = SO(V_{\vp})= SO_2$  which is naturally  extended  to  the homomorphism $\vp : H =  H^s \cdot  T^k   \to  T^1$ , which  sends  the  semisimple part $H^s$  of  $H$   into identity,    and   by identification   of  the tautological $SO_2$-module $V_{\vp} =  \mathbb{R}^2$ with   $\mathbb{C}$ by  choosing one  of  the  two invariant  complex  structures $\pm J$. In  this  case, the  singular orbit  $S_0 = G/H$  is identified  with  the projectivisation
  $P M_{\vp} =  G \times_H P\mathbb{C} = G/K $ of  the   vector  bundle.\\

 Let $\beta_1, \cdots , \beta_p$  be simple  black  roots  (from  $\Pi_B$)  and   $\pi_1, \cdots, \pi_p$ be the associated   fundamental  weights.
 Then  the character  $\chi : Z(H)= T^k \to T^1$   is    determined  by   an infinitesimal   $T$-character
     $\dot{\chi} \in  P_T := \mathrm{span}_{\mathbb{Z}} \Pi_B$  (see  \cite{A-Ch})  and    has   the   form
   $$ \chi(\mathrm{exp}(2 \pi t ) )= \mathrm{exp}(2 \pi \dot{\chi}(t)) ,\,  t \in Z(\mathfrak{h}).  $$
 Sometimes, we  will identify $\chi$  with  $\dot{\chi}$.\\

 \subsubsection{Case of  vector  bundles of  rank  $m >1$}

 The   description  of    standard  vector  bundles  of  rank $m >1$ over  a   flag manifold   $S_0 = G/H$ is  similar  to  the  case of   line  bundles.\\

 Let  $(S_0 =G/H, J^S)$ be  a  flag manifold   associated  with a painted
    Dynkin  diagram  $\Pi = \Pi_B \cup \Pi_W$.  We   fix  a  connected   component  of the white  subdiagram  $\Pi_W$   which is  a string of length $m-1$,  i.e. has  the  type   $ A_{m-1}$  and  corresponds  to a $\mathfrak{su}_m$ ideal of  $\mathfrak{h}$. We  have  the   following  decomposition   of  $\mathfrak{h}$ into a  direct  sum  of ideals
    $$ \mathfrak{h} = \mathfrak{su}_m \oplus \gn' \oplus  \gt^k  $$
where    $\gt^k = Z(\gh)$ is  the  center   and  $\gn'$ is    the  semisimple  ideal complementary  to  $\gsu_m$.\\
 As   in the case  $m=1$,    an   admissible  bundle   is  defined  by  a character  $\chi:  Z(H)= T^k \to  T^1 = e^{i \bR}$    which   determines   the homomorphism
  $\vp : H =  SU_m \cdot N' \cdot T^k  \to  V_{\vp}$   where $V_{\vp} =\bC^m$ is  the tautological  $  SU_m$-module extended   to an $H$-module  by the  conditions
   $$\vp(N') = \mathrm{id}, \ \  \vp |_{T^k}  = \chi$$
     where  $e^{ia} \in \chi(T^k) \, ,  a\in \bR$, acts  on
 $V_{\vp} = \bC^m$  by   complex  multiplication.\\
  Note  that   we  fix  one of  the  two   invariant  complex  structures $J^{V_{\vp}}$ in the $H$-module $V_{\vp}$.  Together  with a complex  structure  $J^S$   on  the base  $S_0$ of the vector  bundle  $M_{\vp} = G \times_H V_{\vp}$, this  defines  a   projectable  invariant  complex  structure  in  the total  space  $M_{\vp}  = G \times_H V_{\vp}$, hence  also  an invariant    complex  structure  $J^F$   on  the  flag manifold  $F = G/K = G/ L \cdot  T^1$  which is  the projectivisation  of  the  vector  bundle  $M_{\vp} = G \times_H V_{\vp}$.\\
Note  that    the  opposite   complex  structure  $- J^{V_{\vp}}$  defines another  projectable   complex  structure  $J'$ on  $M_{\vp}$    and  another  invariant  complex  structure  $(J')^F$ on  $F$.\\

The  following   definition  describes the  data  which  determine  an   admissible  homogeneous  vector   bundle  $M_{\vp} = G\times_H V_{\vp}  \to  S_0 =G/H$  of  rank  $m>1$  together  with  an invariant  complex  structure $J$.

\bd Let   $\Pi =\Pi_B \cup \Pi_W$   be  a painted Dynkin  diagram which  defines  a  flag manifold   $(S_0 =G/H,J^S)$.   A  triple   $(A_{m-1}, \chi, \beta)$  ,  where   $A_{m-1} =\{\a_1,  \cdots, \a_{m-1}\}$ is  a  string, i.e.   a   connected   component of  the   white  subdiagram   $\Pi_W$   of  type  $A_{m-1}$,  $\chi : Z(H) =T^k \to T^1$   a  character    and  $\beta$    is one of  the  end  roots of   $A_{m-1}$, ( the left $\beta= \a_1$  or  the right
 $\beta = \alpha_{m-1}$  )  is  called    a  data.
\ed

\bp \label{ pair} \cite{AZ}
  Let  $(S_0 =G/H, J^S)$ be  a  flag manifold  with  reductive  decomposition   $\gg = \gh + \gm$, associated  with a painted
    Dynkin  diagram  $\Pi = \Pi_B \cup \Pi_W$. A  data  $(A_{m-1}, \chi, \beta)$    defines
    an admissible homogeneous vector bundle $M_{\vp} = G\times_H V_{\vp}$ with  a  complex  structure $J$.
      The  restriction  of the  complex  structure $J$   to  $M_{reg}  =  G/L \times \bR^+$   is  defined as  follows.
      The   $B$-orthogonal  reductive  decomposition  of a regular orbit  $ G/L$    can be  written  as
      $$\gg = \gl + (\bR Z^0_F + \gp) = \gl +  (\bR Z^0_F + \gq  + \gm) $$
       where $ \gh = \gk  + \gq,\ ,   \gk =  \gl + \bR Z^0_F$
       and  $Z^0_F$ is  the fundamental vector of  the  principal  $T^1$ bundle  $G/L \to F =G\times_H PV_{\vp}= G/K $, normalised  by   $B(Z^0_F, Z^0_F)= -1$ .
        The  complex  structure  $J^F$ induces  the  invariant   CR  structure
        $(\mathcal{H}, J^{\mathcal{H}})$ in  $G/L$. It is  extended  to  the invariant  complex  structure $J$ on  $M_{reg}  = G/L \times \bR^+$ by the   formula
        $$  J \partial_t = \frac{1}{a(t)} Z^0_F , \ \  J Z^0_F  = - a(t)  \partial_t$$
\noindent where $a(t)$ is a non-vanishing function.

The  $B$-orthogonal  reductive   decomposition   of the  singular orbit   $S_0 = G/H$  and   a  regular orbit  $G/L$
 can be  rewritten  as

 \begin{equation}\label{standard decomposition of $S_0$}
\gg  = \gh + \gm = (\gsu_{m} \oplus  \gn' \oplus \gt^k) + \gm =
  (\gsu_m \oplus \bR Z^{\chi}  \oplus \gn' \oplus \mathrm{ker}\chi ) + \gm.
\end{equation}

\begin{equation}\label{standard decomposition of $S_t$}
\gg = \gl +  \gp  =  (\gu_{m-1} \oplus  \gn'  \oplus \mathrm{ker}\chi ) +  (\bR Z^0_F + \gq + \gm  )
\end{equation}
where  $Z^{\chi}   = B^{-1} \dot{\chi},\,  \gt^k = Z(\gh)= \bR Z^{\chi} + \mathrm{ker \dot{\chi}},     \gk =  \gl + \bR Z^{\chi}$. We identify   the    subalgebra  $\gu_m =\gsu_m + \bR Z^{\chi}$  with $\vp(\gu_m) = \gu(V_{\vp})$  and denote  by  $\gu_{m-1}$  the  stabilizer of  a  fixed vector $e_0 \in  V_{\vp} = \bC^m$; Finally, $\gq$ is   the  invariant  complement  to  $\gu_{m-1}$ in  $\gu_m = \gu_{m-1}  + \gq$.
\ep

  \subsection{Invariant K\"ahler  structures on the  total  space  $M_{\vp}$ of  a     standard  vector  bundle}

    Invariant K\"ahler structures on  a standard \C1 manifold $M_{\vp} = G \times_H V_{\vp}$ are described by segments (an interval  or  a  ray) in the $T$-Weyl chamber
    $$C(J^F) = \{\beta>0, \beta_1 >0, \cdots , \beta_k >0\} \subset i Z(\mathfrak{k})=\mathfrak{t}^k + \mathbb{R}Z^0$$
     of  the  flag manifold $F = G/K$ corresponding to the complex structure $J^F$ (see Theorem 7 in the Appendix). Here
\be\label{Z^0def}
Z^0 = -i Z^0_F
\ee
where   $Z^0_F$  is  the  fundamental vector.\\
 We may  assume  that   $\beta(Z^0) >0$ where
$\beta$ is the new black root in the Dynkin diagram of $G/K$.\\
Choose   a vector $ Z_0 \in i Z(\mathfrak{k}) $ such  that $\beta(Z_0)=0,\,  \beta_i(Z_0)>0,\,  i =1, \dots, k$. Geometrically,  the vector  $Z_0$  belongs  to  the  face   $\beta =0$   of  the  Weyl  chamber $ C(J^F)$  and  its  projection  to  $ i Z(\mathfrak{h})$  is in  the  Weyl  chamber  $C(J^S)$.\\

\bd  A  segment  (an interval  or  a ray) in  $C(J^F)$ of   the  form  $(Z_0 Z_d), \, \beta(Z_0) =0 $,  which is  parallel  to  the  fundamental vector  $Z^0$   together with   a  parametrization   $Z_0 + f(t) Z^0$  such  that   $\dot{f}(t) >0, f(0)=0,  Z_d = Z_0 + f(d) Z^0$, is  called   an {\it admissible  segment}.\\
 \ed

 \bt   (\cite{AZ}, Proposition  17 ,  see  also \cite{P-S})
Let  $(M_{\varphi}, J)$  be an admissible  vector  bundle  associated  with  a data  $(A_{m-1},\chi,\beta)$. Any   admissible segment $(Z_0 Z_d)  \subset C(J^F) $
     defines  a K\"ahler  metric in the  tubular $S_0$-punctured neighbourhood
    $M   =  (0,d) \times G/L \subset  M_{\vp} \setminus S_0$   of  the   zero  section $S_0$  of  the  vector  bundle $M_{\vp} \to S_0 = G/H$ given  by

$$ g_{reg} = dt^2 + (\dot f \theta^0)^2 + \pi_F^* g_0 + f(t) \pi_F^*g^0.$$

    Here $\pi_F : M_{reg}=  G/L \times \bR^+ \to  F  =G/K$ is   the  natural projection and  $g_0 =  -\omega_{Z_0}\circ J^F $, $g^0 = -\omega_{Z^0}\circ J^F$, where $\omega_{Z_0}$, $\omega_{Z^0}$ are the closed invariant forms on $F$  associated  with   $Z_0, Z^0$ (see Proposition \ref{propcorrespmetrics} in the Appendix below for the correspondence between vectors in $C(J^F)$ and forms on $F$). Any  invariant  K\"ahler metric of  standard  type  can be obtained  by  this  construction.\\
   The K\"ahler metric  $g$   smoothly   extends  to   the  zero   section $S_0$ if and  only if   the function $f(t)$
   is  extended  to a  smooth even   function on  $\mathbb{R}$ which satisfies  the   following  Verdiani  conditions  \cite{V} :
$$f(0) = \dot f(0)=0,\,  \ddot f(0) =\kappa,$$
where
\begin{equation}\label{DefKappa}
\kappa = 2 \pi/T_0, \ \ \ T_0 = \min \{t > 0 \ | \ \exp(t Z_0) \in L \}
\end{equation}

 Moreover, the  invariant K\"ahler metric   $g$ is  geodesically   complete  on  $M_{\vp} = G \times_H V_{\vp}$ if  and only if  the  function $f(t)$ is defined on $\mathbb{R}^+$  and   satisfies  the  Verdiani  conditions.

\et
\medskip

 Finally, we recall the conditions  for the  K\"ahler metric    associated   with  an  admissible  segment $(Z_0 Z_d)$  to be a  K\"ahler-Einstein  metric.

\bt \label{maintheoronKE}(Theorem 34 in \cite{AZ})
Let $M_{\vp}$ be a standard \C1 manifold, i.e. the total space of an admissible bundle $M_{\vp} = G \times_H V_{\vp} \rightarrow S_0$ over the singular orbit $(S_0 = G/H, J^S)$  and
 $(F = G \times_H P V_{\vp} = G/K, J^F)$ be  the flag manifold associated with regular orbits.
 The invariant  K\"ahler metric  $g$ in $M_{\vp}$   associated  with  an  admissible  segment $(Z_0 Z_d) \subseteq C = C(J^F)$ in   the  $T$-Weyl  chamber  $C(J^F)$ is  a K\"ahler-Einstein metric  with Einstein  constant  $\lambda$ if  and only if \\
 (i) the Koszul vector
 $Z^{Kos} \in C(J^F)$ (which defines  the  invariant K\"ahler-Einstein metric  on the  flag manifold  $(F,J^F)$, see the Appendix),  the  initial vector  $Z_0$ of  the  segment   and  the  fundamental vector $ Z^0$   are  related  by
  \be\label{algCond}
Z^{Kos} = \lambda Z_0 + \kappa m Z^0
\ee
\noindent where $m = \dim(V_{\vp})$ and  $\kappa$ is defined by (\ref{DefKappa});
 \\
  (ii)    the  function  $f(t)$   satisfies   the   equation

   \be \label{EinstEq}
   \ddot f(t) + \frac12 A(f) \dot f^2 + \lambda f =\kappa m
   \ee

\noindent with  the initial  conditions\\
   $\lim_{t \to 0}f(t) = \lim_{t \to 0}\dot f(t)=0,\, \lim_{t \to 0}\ddot f(t)= \kappa, $

\noindent where $A(f) = \sum_{\alpha \in R_{\mathfrak{m}}^+} \frac{\alpha(Z^0)}{\alpha(Z_0) +f \alpha(Z^0)}$  and $R_{\mathfrak{m}}^+$  is  the set of the positive black roots of $G/K$,  see   Appendix.\\
Moreover, the K\"ahler-Einstein metric can be extended to a complete metric if and only if $\lambda \leq 0$, and  the segment extends to a ray $Z_0 + \bR^+ Z^0$ in $C(J^F)$.
\et

\noindent The proof,  given  in  \cite{AZ}, is  based on the following

\bt \label{theoremDescribesKE}
If the condition (\ref{algCond}) of Theorem \ref{maintheoronKE} is fulfilled, then the function $f(t)$ parametrizing the segment $(Z_0 Z_d)$ which gives the K\"ahler-Einstein metric is the inverse to the function

\be\label{t(f)}
t(f) = \int_{0}^{f} \sqrt{\frac{P(s)}{2 \int_0^s (\kappa m - \lambda v) P(v) dv}} ds
\ee

\noindent where $P$ is the polynomial defined by $P(x) = \Pi_{\alpha \in R_{\mathfrak{m}}^+} (\alpha(Z_0) + x \ \alpha(Z^0))$.
\et
  {\bf  Remark}  If   the  necessary  and  sufficient  conditions    are   fulfilled this  theorem reduces  the  explicit   construction  of a K\"ahler-Einstein metric  to the construction  of  the  inverse function  $f(t)$   to  the   function  $t(f)$.

\medskip

\subsection{The  main  results}

 Let  $(F = G/K, J^F)$  be  the  flag manifold   with  an invariant  complex  structure associated  with  a painted  Dynkin  diagram   $\Pi = \Pi^F_B \cup \Pi^F_W$.  Denote  by  $\{  \beta_0, \beta_1, \cdots , \beta_p \} = \Pi^F_B$  the  simple  black  roots   and by
 $\pi_0, \pi_1, \cdots , \pi_p$ the  associated     black fundamental  weights. The  Koszul  form  $\sigma_F = B \circ Z^{Kos}$  associated   with  the  Koszul  vector $Z^{Kos}$
   admits  a  decomposition

 \be \label{ Koszul numbers}
B\circ Z^{Kos} = n_0 \pi_0 + n_1 \pi_1 + \cdots  + n_p \pi_p
 \ee
 where the natural numbers   $n_i$  are  called   the  {\it Koszul  numbers  of   the complex   flag manifold  $(F, J^F)$}.\\

   Now we are ready to state our main  theorems  which  give  necessary  and  sufficient  conditions in order  that the  admissible  vector  bundle  $M_{\vp} = G \times_H V_{\vp}  \to  S_0 = G/H$   over   a  flag  manifold $(S_0 =G/H, J^S)$ associated with a painted  Dynkin  diagram  $\Pi = \Pi_B \cup \Pi_W$  admits  an  invariant
   K\"ahler-Einstein  metric,  where  $G$ is   one  of  the  classical  compact  Lie  groups   $SU_n, Sp_n, SO_n$.  Recall   that a \C1  $G$-manifold  $M_{\vp}$ having $G/H$ as singular orbit and endowed with  a  complex  structure   $J$ is  defined   by  the  data
   $(\mathfrak{su}_{m}, \chi, \beta)$,   where $\mathfrak{su}_{m}$ is  a   connected   component  of  the   white   subdiagram $\Pi_W$ of $G/H$, $\beta = \beta_0$  is one of  the  end roots of  the   string $\mathfrak{su_{m}}$  and $\chi : Z(H) = T^k \to T^1$ is  a  character.\\
  These  data   define  a  complex  structure   on  the  flag  manifold   $F = G \times_H PV_{\vp} = G/K  $
such  that   $(F, J^F)$ corresponds  to  the painted Dynkin diagram  $\Pi = \Pi^F_B \cup \Pi^F_W$  obtained   form  the painted Dynkin  diagram  $\Pi = \Pi_B \cup \Pi_W$ of  $S_0 = G/H$  by painting  the   simple  root $\beta= \beta_0$ into   black. So  $\Pi^F_B  = \Pi_B \cup \{\beta_0 \}$.\\

For   classical simple  Lie  algebras $\mathfrak{g}$ of  types  $A, B, C, D $ we  use  the    standard notation  for the root  system  $R_{\mathfrak{g}}$  and   the simple  root  system   $\Pi_{\mathfrak{g}}$  as in   \cite{G-O-V}:
$$
\begin{array}{llll}
  R_A =& \{\e_i -\e_j\},\, & R_B = & \{ \e_i -\e_j,\pm \e_i\},\, \\
    \Pi_A = & \{ \a_i :=\e_i - \e_{i+1} \} & \Pi_B = & \{\a_i:= \e_i - \e_{i+1}; \e_{\ell} \} \\
    R_C = &\{ \e_i -\e_j,  \pm 2 \e_i\} \, & R_D =& \{\pm \e_i \pm \e_j\} \\
  \Pi_C = &\{ \a_i:=\e_i - \e_{i+1}; 2\e_{\ell}\} &  \Pi_D = &\{ \a_i:= \e_i - \e_{i+1}; \a_{\ell}:= \e_{\ell-1}+ \e_{\ell} \}
\end{array}
  $$
where  $\ell$ is  the  rank.
Now  we  are  ready to  state  the main  results of  the paper.
The  following   Theorem \ref{prop triple2} (resp. Theorem \ref{prop triple4})  gives  necessary  and  sufficient  conditions  for  an  admissible  vector  bundle  $M_{\vp} = G \times_H V_{\vp} $  of  rank   $m = \dim  V_{\vp} >1$ (resp. $m=1$) to  admit  a K\"ahler-Einstein standard invariant metric.

\bt \label{prop triple2}
  Let   $(S_0 = G/H, J^S)$  be   the  flag manifold  of one of  the   classical simply   connected  Lie groups    $ G= SU_n, Sp_n,\widetilde{SO}_n = Spin_n $     defined  by  the painted  Dynkin  diagram $\Pi = \Pi_W \cup \Pi_B, \, \Pi_B = \{ \beta_1, \cdots, \beta_p \}$, and let $n_1, \dots, n_p$ be the Koszul numbers of $G/H$.

Let  $m>1$ an integer, let $A_{m-1} = \{\a_1, \cdots, \a_{m-1}\}$ be a  white  string of  $\Pi_W$ and $\beta = \a_1$ (resp. $\beta = \a_{m-1}$).  Then, we have the following

\begin{enumerate}
\item[(i)]  the admissible vector bundle $M_{\vp} =  G \times_H V_{\vp} \to  S_0$  of  rank  $m$  associated  with  the data $(A_{m-1}, \chi, \beta)$ admits a Ricci- flat  K\"ahler  standard  metric if and only if  $m|n_j$ for $j=1, \dots, p$  and

\begin{equation}\label{chienunciatoTh11}
\chi = \sum_{j = 1}^p \frac{n_j}{m} \pi_j  \ \ {\rm (resp. \ \ } \chi = - \sum_{j = 1}^p \frac{n_j}{m} \pi_j)
\end{equation}

\item[(ii)] the admissible vector bundle $M_{\vp} =  G \times_H V_{\vp} \to  S_0$  of  rank  $m$  associated  with  the  triple $(A_{m-1}, \chi, \beta)$, where $\chi = k_1 \pi_1 + \cdots + k_p \pi_p$ and $\beta = \alpha_1$ (resp. $\beta = \a_{m-1}$) admits a unique K\"ahler-Einstein standard   metric $g$ defined  in  a neighborhood  of  the   singular   section with Einstein constant $\lambda > 0$ if and only if $k_j < \frac{n_j}{m}$ (resp. $k_j > -\frac{n_j}{m}$) and with Einstein constant $\lambda < 0$ if and only if $k_j > \frac{n_j}{m}$ (resp. $k_j < -\frac{n_j}{m}$).

In the case $\lambda < 0$ the  metric is  extended  to a   globally  defined  complete metric  in  $M_{\vp}$.

\end{enumerate}

\et

\bt \label{prop triple4}
  Let $G$ and $S_0$ be as in Theorem \ref{prop triple2}. Then the admissible vector bundle $M_{\vp} =  G \times_H V_{\vp} \to  S_0$ of  rank  $m = 1$ with $S_0$ as only singular orbit associated to the infinitesimal character $\dot{\chi}=   k_1 \pi_1 + \cdots + k_p \pi_p$ admits a K\"ahler-Einstein standard \C1 structure with Einstein constant $\lambda =0$ (resp. $\lambda > 0$, $\lambda < 0$) if and only if $k_j = n_j$ (resp. $k_j > n_j$, $k_j < n_j$), where the $n_j$'s, $j=1, \dots, p$, denote as above the Koszul numbers of $G/H$.
\et

 Note  that the last theorem includes the case when $S_0$ is the manifold of  full   flags (i.e. $H = T^{\ell}$ is  a maximal  torus).

\subsection{Calculation of Koszul numbers  and  examples }

\noindent Let us recall the flag manifolds $F = G/K$ of the classical groups $G$
(see, for example, \cite{A-P}, \cite{arv}):

\begin{enumerate}

\item[-] $SU(n)/S(U(n_1) \times \cdots \times U(n_s) \times U(1)^m)$

$n = n_1 + \cdots + n_s + m$, $s, m \geq 0$

\item[-] $SO(2n+1)/U(n_1) \times \cdots \times U(n_s) \times SO(2r+1) \times U(1)^m$

\item[-] $Sp(n)/U(n_1) \times \cdots \times U(n_s) \times Sp(r) \times U(1)^m$

\item[-] $SO(2n)/U(n_1) \times \cdots \times U(n_s) \times SO(2 r) \times U(1)^m$

$n = n_1 + \cdots + n_s + m + r$, $s, m, r\geq 0, r \neq 1$

\end{enumerate}

\noindent The Koszul numbers $n_j$   for $F$ endowed with a $G$-invariant complex structure $J^F$ are determined    by   the     corresponding  painted Dynkin  diagram  as  follows (\cite{A-P}): \\
$n_j = b_j+2$, where $b_j$ equals the number of white roots  connected to the black root $\beta_j$ ,  with  the  following  exceptions.\\
  For  the  group $G = SO_{2n+1}$ of  type  $B_n$,   each  long  root   of  the  last  white  chain
    which defines the root system of the  type $so_{2r+1}$  is counted as two.\\
   For $G$  of  type  $C_n$,    each  root  of  the last white  chain    of type   $sp_r$  is  counted  as  2.\\
   For  $G$ of  type $D_n$,   the last white chain  of  type  $ so_{2r}$ is
considered as a chain of length 2(r-1).\\
 If r = 0 and one of the two end  roots is
white and the other one is black,  the Koszul number associated to the  end     black   root $\beta$ is 2(k-1), where k is the number of white roots connected with   $\beta$.  \\

\smallskip
\begin{Ex}\label{esempipiattoealtri} \rm
\noindent Let us consider for example  the  flag manifold  $(G/H,J)$ given by the following painted diagram:

\begin{equation}\label{paintedG/H}
\circ - \circ - \bullet - \circ - \circ - \bullet - \underset{}\circ - \underset{}\circ - \underset{}\circ - \underset{}\circ - \circ
\end{equation}

\noindent The Koszul numbers associated to the black roots are

$$n_1 = 6, \ n_2 = 9$$

For the first white string on the left $\mathfrak{su}_{m}$, $m=3$, condition $m|n_j$ in Theorem \ref{prop triple2} is satisfied, so there exists a Kahler-Einstein admissible vector bundle of rank $3$ with Einstein constant $\lambda = 0$ if we choose data $(\mathfrak{su}_{m}, \chi, \beta)$ with $\mathfrak{su}_m$ being this string (both when the new black root $\beta$ is the first and the second node of the string).
For the white string $\mathfrak{su}_{m}$, $m=6$ on the right, the condition $m|n_j$ in Theorem \ref{prop triple2} is not satisfied since $6$ does not divide $9$, so the admissible vector bundle corresponding to the choice of this string does not admit a Ricci-flat structure.

\medskip

\noindent If $G/H$ is given by the following painted diagram:

\begin{equation}\label{paintedG/H}
\bullet - \circ - \circ - \circ - \bullet
\end{equation}

\noindent the Koszul numbers associated to the black roots are

$$n_1 = n_2 = 5$$

For the central white string $\mathfrak{su}_{m}$, $m=4$, condition $m|n_j$ in Theorem \ref{prop triple2} is not satisfied, so the admissible vector bundle corresponding to the choice of this string does not admit a Ricci-flat structure.

\end{Ex}

\begin{Ex}
\rm The conditions for the existence of a Kahler-Einstein metric (with Einstein constant $\lambda=0, \lambda>0$ or $\lambda <0$) are satisfied in particular when the painted Dynkin diagram of the singular orbit $S_0 = G/H$ consists of a white $A_{m-1}$ string only (i.e. $S_0 = SU(m)/SU(m)$ is a point), since in this case the Koszul numbers $n_j$, $j=1, \dots, p$ of $G/H$ all vanish. So there exists a Kahler-Einstein  standard cohomogeneity one $SU(m)$-manifold having a point as singular orbit for any value of the Einstein constant $\lambda$.

\noindent In order to determine explicitly the metric for any case, assume for example that $\beta= \a_1 =  \varepsilon_1 - \varepsilon_2$ (the case when $\beta = \a_{m-1} = \varepsilon_{m-1} - \varepsilon_m$ is similar) and observe that $Z(\mathfrak{h}) = \{ \bar 0 \}$ and that the set of positive black roots of $G/K$ is

$$R_{\mathfrak{m}}^+ = \{ \beta, \beta + (\varepsilon_2 - \varepsilon_3), \dots, \beta +\cdots +  (\varepsilon_{m-1} - \varepsilon_m) \}$$

\noindent Then we have that the polynomial $P(x)$ is

$$P(x) = \Pi_{\alpha \in R_{\mathfrak{m}}^+} (\alpha(Z_0) + x \ \alpha(Z^0)) = \Pi_{\alpha \in R_{\mathfrak{m}}^+} x \ \alpha(Z^0) = \beta^0 x^{m-1}$$

\noindent (where we are denoting $\beta^0 = \beta(Z^0)^{m-1}$) and then

$$\int_0^s (km - \lambda v) P(v) dv = \beta^0 k s^m - \lambda \frac{\beta^0 }{m+1} s^{m+1}.$$

\noindent So, by Theorem \ref{theoremDescribesKE}, the function $f(t)$ which determines the metric is the inverse to

$$t(f) = \int_{0}^{f} \sqrt{\frac{P(s)}{2 \int_0^s (\kappa m - \lambda v) P(v) dv}} ds  = \int_{0}^{f} \sqrt{\frac{1}{2 \kappa s - \lambda \frac{2 }{m+1} s^2}} ds$$

\noindent By a straight calculation, one then sees

\begin{enumerate}
\item[- \ $\lambda=0$:] $$t(f)= \int_{0}^{f} \sqrt{\frac{1}{2 \kappa s }} ds = \frac{\sqrt2}{\kappa} \sqrt f$$
\noindent so that $t \in [0, + \infty)$ and
$$f(t) = \frac{k^2}{2} t^2$$
\noindent This is the flat metric on $\bC^m$, endowed with the canonical $SU(m)$ action.
\item[- \ $\lambda>0$:] $$t(f)= \int_{0}^{f} \sqrt{\frac{1}{2 \kappa s - \lambda \frac{2 }{m+1} s^2}} ds = -\frac{2}{\sqrt{-b}} artg \frac{1}{\sqrt{-b}} \sqrt{\frac{bf}{f + \frac{a}{b}}}$$
\noindent so that $t \in [0, \frac{\pi}{\sqrt{-b}}]$ and
$$f(t) = -\frac{a}{b}\sin^2(\frac{\sqrt{-b}}{2}t)$$
\noindent where we are denoting $a = 2 \kappa$, $b= - \lambda \frac{2 }{m+1}$: this is the (non complete) Fubini-Study metric on $\bC^m$.
\item[- \ $\lambda<0$:] $$t(f)= \int_{0}^{f} \sqrt{\frac{1}{2 \kappa s - \lambda \frac{2 }{m+1} s^2}} ds = \frac{1}{\sqrt{b}} ln \frac{1+ \sqrt{\frac{f}{f + \frac{a}{b}}}}{1 - \sqrt{\frac{f}{f + \frac{a}{b}}}}$$
\noindent so that $t \in [0, + \infty)$ and
$$f(t) = \frac{a}{b}\sinh^2(\frac{\sqrt{b}}{2}t)$$
\noindent where we are denoting $a = 2 \kappa$, $b= - \lambda \frac{2 }{m+1}$: this is the hyperbolic metric on the open disk endowed with the canonical action of $SU(m)$.
\end{enumerate}

\end{Ex}

\section{Proofs}

\noindent The proofs of Theorem \ref{prop triple2} and Theorem \ref{prop triple4} consist in finding the conditions under which there exists a Lie algebra character $\chi: Z(\mathfrak{h}) \rightarrow \bC$ such that the above algebraic condition $Z^{Kos} = \lambda Z_0 + \kappa m Z^0$ (\ref{algCond}) in Theorem \ref{maintheoronKE} is satisfied. In order to do this, we need to calculate $Z^0$ and $Z^{Kos}$.

\bl\label{LemmaKZ^0}

Let $G$ be a simply connected group with Lie algebra $\mathfrak{su}_n$, $\mathfrak{sp}_n$, $\mathfrak{so}_{2n}$, $\mathfrak{so}_{2n+1}$ and let $S_0 = G/H$ be a flag manifold with painted Dynkin diagram $\Pi = \Pi_B^H \cup \Pi_W^H$.

\noindent Let $G \times_H V$ be the standard admissible bundle of rank $m>1$ defined by the data $(A_{m-1}, \chi, \beta)$, where $A_{m-1} =\{\a_1, \dots, \a_{m-1} \}$ is a white string in $\Pi_W^H$,  $\chi: Z(\mathfrak{h}) \rightarrow \mathbb{C}$ is a Lie algebra character and $\beta = \alpha_1$ (resp. $\beta = \alpha_{m-1})$ is the new black root in the painted Dynkin diagram of the flag $G/K$ associated to the regular orbits.

\noindent Let $Z^0, \kappa$ be defined by (\ref{Z^0def}) and (\ref{DefKappa}). If $\pi_{s}$, $\pi_{s+1}$ denote the fundamental weights of the black roots $\beta_s$, $\beta_{s+1}$ of the diagram of $G/H$ connected to $\alpha_1$ and $\a_{m-1}$ respectively, then, up to sign, we have

$$\kappa Z^0 = B^{-1}(\chi + \pi_0 - \frac{m-1}{m} \pi_{s} - \frac{1}{m} \pi_{s+1})$$
 \begin{equation}\label{KZ^0INLEMMA22}
 {\rm (resp.} \ \  \kappa Z^0 = B^{-1}(\pi_0 - \chi - \frac{m-1}{m} \pi_{s+1} - \frac{1}{m} \pi_{s}))
\end{equation}
\noindent with the exception of the following two cases:

\begin{enumerate}

\item[(1)] $\mathfrak{g} = \mathfrak{so}_{2n+1}$ and the painted Dynkin diagram of $G/H$ is

\begin{equation}\label{paintedG/H1lemma}
\cdots -  \underset{\beta_s}\bullet - \underset{\alpha_1}\circ - \cdots - \underset{\alpha_{m-1}}\circ \Rightarrow \underset{\beta_{s+1}}\bullet
\end{equation}

\noindent then

$$\kappa Z^0 = B^{-1}(\chi + \pi_0 - \frac{m-1}{m} \pi_{s} - \frac{2}{m} \pi_{s+1})$$
 \begin{equation}\label{KZ^0INLEMMA24lemma}
 {\rm (resp.} \ \  \kappa Z^0 = B^{-1}(\pi_0 - \chi - \frac{2(m-1)}{m} \pi_{s+1} - \frac{1}{m} \pi_{s}))
\end{equation}

\item[(2)] $\mathfrak{g} = \mathfrak{so}_{2n}$ and the painted Dynkin diagram of $G/H$ is

\begin{align}\label{paintedG/H2lemma}
 \cdots - \underset{{\beta_s}}{\bullet} - \underset{{\alpha_1}}{\circ} - \dotsb - \underset{{}}{\overset{\overset{\textstyle\bullet_{{\beta{s+1}}}}{\textstyle\vert}}{\circ}} \,-\, \underset{{\alpha_{m-1}}}{\circ} \\
\end{align}

\noindent then

$$\kappa Z^0 = B^{-1}(\chi + \pi_0 - \frac{m-1}{m} \pi_{s} - \frac{2}{m} \pi_{s+1})$$
 \begin{equation}\label{KZ^0INLEMMA23lemma}
 {\rm (resp.} \ \  \kappa Z^0 = B^{-1}(\pi_0 - \chi - \frac{m-2}{m} \pi_{s+1} - \frac{1}{m} \pi_{s}))
\end{equation}

\end{enumerate}

\noindent where we are denoting $B^{-1}(\xi)$ the dual of $\xi$ with respect to the Killing form $B$, that is $\xi(X) := B(B^{-1}(\xi), X)$.

\smallskip

\noindent For the admissible vector bundle $M_{\vp} =  G \times_H V_{\vp} \to  S_0$ of  rank  $m = 1$ with $S_0$ as only singular orbit (i.e. $G/K = G/H$) defined by the pair  $(A_{m-1}, \chi)$ we have

\begin{equation}\label{KZ^0INLEMMA1}
\kappa Z^0 = B^{-1}(\chi)
\end{equation}

\el

\begin{remark}
\rm \noindent If $\alpha_1$ is the first (resp. $\alpha_{m-1}$ is the last) node of the diagram, then we have no black root $\beta_s$ (resp. $\beta_{s+1}$) and in formulas (\ref{KZ^0INLEMMA22}), (\ref{KZ^0INLEMMA24lemma}) and (\ref{KZ^0INLEMMA23lemma})  the term in $\pi_s$ (resp. in $\pi_{s+1}$) cancels.
\end{remark}

\noindent Before starting to prove the Lemma, let us fix some notation which will be fundamental in the proof.

\noindent As we have recalled above, the  stability  subgroup $L$ of a regular  orbit in
$M_{\vp} = G \times_H V_{\vp}$ can be   identified  with  the   stability  subgroup $H_e$
of  a non-zero  vector $e \in V_{\vp}$  and  the  corresponding
stability   subgroup $K$ of  the  flag manifold $F=G \times_H
PV_{\vp}$  with the    stabilizer $H_{[e]}$ of the line $[e] \in
PV_{\vp}$.

\noindent This holds true also when $m= \dim(V_{\vp}) =1$, in which case $G/K = G/H$, that is we have no new black root $\beta$.

\smallskip

Given the Lie algebra character $\chi: Z(\mathfrak{h}) \rightarrow \mathbb{C}$, let us denote $\ga = \mathrm{ker}\, \chi$ and define the
following direct sum orthogonal decomposition

\be \label{decompositionof"h"}
 \mathfrak{h}=\mathfrak{su}_m + \mathfrak{n}' + Z(\mathfrak{h})=
 \mathfrak{su}_m + \mathfrak{n}' + \mathfrak{a} + {\mathbb{R}}Z^{\chi} \ee
 where $Z^{\chi}$  is   the  vector in  $\ga^{\perp}$ such
 that  $\chi(\exp t Z^{\chi}) = e^{it}$.

 We  identify $V_{\vp}$  with the  Hermitian  space $\mathbb{C}^m$ such that   the
 standard basis $e_j, \, j=1,\dots, m$   consists  of  weight vectors
 with weights $\varepsilon_j$   w.r.t. the Cartan subalgebra  $\vp(\gc)$
 and   the  simple  roots  $\alpha_j \in A_{m-1}$ satisfy  $\alpha_j|_{\mathfrak{c}} =\varepsilon_i -
\varepsilon_{i+1}$.

\noindent In the case $m > 1$, we have either $\beta=\alpha_1$, in which case we take $e = e_1$ and
set
\begin{equation}\label{Zbeta1}
Z^{\beta}  := i \mathrm{diag}((m-1), -\mathrm{id}_{m-1}) \in  \mathfrak{su}_m \subset \mathfrak{h}
\end{equation}

\noindent or $\beta = \alpha_{m-1}$, in which case we choose $e = e_m$ and set
\begin{equation}\label{Zbeta2}
Z^{\beta}  := i \mathrm{diag}(-\mathrm{id}_{m-1}, m-1) \in  \mathfrak{su}_m \subset \mathfrak{h}
\end{equation}

\noindent (when there is no risk of confusion, with a slight abuse of notation in the following we will denote by $Z^{\beta}$ both the element of $\mathfrak{su}_m$ and its immersion in $\mathfrak{h}$, see also (\ref{Zbeta1imm})-(\ref{Zbeta32}) below).

Since $Z^{\chi}$ goes under the Lie algebra representation to $i \mathrm{id}_{m}$, in both cases the
element
$$Z^{\mathfrak{l}}= Z^{\beta} -(m-1)Z^{\chi}$$
 annihilates $e$, hence it
 belongs to  the  stability   subalgebra $\mathfrak{l}$. Then the fundamental vector $Z^0_F$ coincides with the vector of  the plane
$\mathrm{span}(Z^{\beta}, Z^{\chi})$ orthogonal to
$Z^{\mathfrak{l}}$ and normalized by $B(Z^0_F, Z^0_F)=-1$.
Recall that the relations between the stability  subalgebras
  $ \mathfrak{h}, \mathfrak{k},\, \mathfrak{l}$  of the flag manifolds  $S_0 =G/H, F= G/K$, the CR manifold
  $G/L$   and  their centers are given by

$$
\begin{array}{lll}
\mathfrak{h} =& \mathfrak{su}_m + \mathfrak{n}' +Z(\mathfrak{h}),&
Z(\mathfrak{h}) = \mathfrak{a} + \mathbb{R}Z^{\chi} \\
$$ \mathfrak{l} = & \mathfrak{su}_{m-1}  + \mathfrak{n}' + Z(\mathfrak{l}),
 & Z(\mathfrak{l})=  \mathfrak{a}+ \mathbb{R}Z^{\mathfrak{l}} \\
\mathfrak{k} =&\mathfrak{su}_{m-1}+  \mathfrak{n}'+Z(\mathfrak{k})=
\mathfrak{l}+ \mathbb{R}Z^0  &
 Z(\mathfrak{k})=  Z(\mathfrak{l})+ \mathbb{R}Z^0  =Z(\mathfrak{h})+\mathbb{R}Z^0
\end{array}
  $$
	
\noindent where we  denote  by $\mathfrak{su_{m-1}}$
the stability   subalgebra of the  vector $e$ in
$\mathfrak{su}_{m}$.

\smallskip

\noindent In the case $m=1$, where as we have observed above $G/H = G/K$ we have no new black root $\beta$, we have $Z(\mathfrak{h}) = Z(\mathfrak{k})$, $Z^{\mathfrak{l}} = \bar 0$ and $Z(\mathfrak{l})=  \mathfrak{a} = \ker \chi$.

\bigskip

{\bf Proof of Lemma \ref{LemmaKZ^0}:} As we have seen above, the vector $Z^0_F$ is given by the $B$-orthogonal decomposition $Z(\mathfrak{k}) = Z(\mathfrak{l}) + \bR Z^0_F$, and $B(Z^0_F, Z^0_F) = -1$, being $Z(\mathfrak{l}) = ker(\chi)$ in the case $m=1$, while for $m > 1$

$$Z(\mathfrak{l}) = \ker(\chi) + \bR Z^{\mathfrak{l}}, \ \ Z^{\mathfrak{l}} = Z^{\beta} - (m-1) Z^{\chi} $$

\noindent where $Z^{\chi} \in Z(\mathfrak{h})$ is orthogonal to $ker(\chi)$ and $\chi(Z^{\chi}) = i$.

\smallskip

\noindent Let $m>1$: by the well-known structure of classical semi-simple Lie algebras, in the case $\mathfrak{g} = \mathfrak{su}_n$, if $\beta = \alpha_1$ (resp. $\beta = \alpha_2$), then by (\ref{Zbeta1}) (resp. (\ref{Zbeta2})) above we have

\begin{equation}\label{Zbeta1imm}
Z^{\beta} = D_m := i \ diag(O, m-1, -id_{m-1}, O)
\end{equation}

\begin{equation}\label{Zbeta2imm}
(\ {\rm resp.} \ Z^{\beta} = D_m := i \ diag(O, -id_{m-1}, m-1,  O) \ )
\end{equation}

\noindent where the order of the zero matrices $O$ depends on the position of the $A_{m-1}$ component in the Dynkin diagram, while for the other classical Lie algebras we have

\begin{equation}\label{Zbeta22}
Z^{\beta} = \left( \begin{array}{cc}
D_m & 0 \\
0 & -D_m
\end{array} \right) \ \rm{for} \ \ \mathfrak{g} = \mathfrak{sp}_{2n}, \mathfrak{so}_{2n}
\end{equation}

\begin{equation}\label{Zbeta32}
Z^{\beta} = \left( \begin{array}{ccc}
D_m & 0 & 0 \\
0 & -D_m & 0 \\
0 & 0 & 0
\end{array} \right) \ \rm{for} \ \ \mathfrak{g} = \mathfrak{so}_{2n+1}
\end{equation}

\noindent where $D_m$ is given either by (\ref{Zbeta1imm}) or (\ref{Zbeta2imm}) depending on the choice of $\beta$.

\noindent We are going to show that

\begin{equation}\label{explicitZ^0}
Z^0_F = \frac{1}{\sqrt{-\frac{1}{\| Z^{\chi}\|^2} - \frac{(m-1)^2}{\| Z^{\beta}\|^2}}} \left(  \frac{Z^{\chi}}{\| Z^{\chi}\|^2} + (m-1) \frac{Z^{\beta}}{\| Z^{\beta}\|^2} \right)
\end{equation}

\noindent where we are using the notation $\|Z\|^2 = B(Z,Z)$. Indeed, $Z^{\chi}$ and $Z^{\beta}$ are orthogonal since $Z^{\chi}$ belongs to $Z(\mathfrak{h})$ which consists of matrices of the kind
\begin{equation}\label{Zhsln}
X_m := i \ diag(O, \theta \ id_m, O)
\end{equation}

 for $\mathfrak{g} =  \mathfrak{su}_n$ and

 \begin{equation}\label{Zhsp}
\left( \begin{array}{cc}
X_m & 0 \\
0 & -X_m
\end{array} \right) \ \rm{for} \ \ \mathfrak{g} = \mathfrak{sp}_{n}, \mathfrak{so}_{2n}
\end{equation}

\begin{equation}\label{Zhso2n+1}
 \left( \begin{array}{ccc}
X_m & 0 & 0 \\
0 & -X_m & 0 \\
0 & 0 & 0
\end{array} \right) \ \rm{for} \ \ \mathfrak{g} = \mathfrak{so}_{2n+1}
\end{equation}

\noindent so the claim is true by (\ref{Zbeta1imm})-(\ref{Zbeta32}) and by recalling that the Killing form $B$ is given by $B(X, Y) = 2n \cdot tr(XY), 2(n+1) \cdot tr(XY), 2(n-1) \cdot tr(XY), (2n-1) \cdot tr(XY)$ for $\mathfrak{g} =  \mathfrak{su}_n,  \mathfrak{sp}_{n},  \mathfrak{so}_{2n},  \mathfrak{so}_{2n+1}$ respectively.

\smallskip

\noindent So we have

$$\left\langle \frac{Z^{\chi}}{\| Z^{\chi}\|^2} + (m-1) \frac{Z^{\beta}}{\| Z^{\beta}\|^2}, \ \ Z^{\beta} - (m-1) Z^{\chi} \right\rangle = - (m-1) + (m-1) = 0$$

\noindent which shows that the vector given by (\ref{explicitZ^0}) is orthogonal to $Z^{\mathfrak{l}}$.

\noindent Moreover, (\ref{explicitZ^0}) is orthogonal to $ker(\chi)$ since both $Z^{\beta}$ and $Z^{\chi}$ are ($Z^{\chi}$ by definition, $Z^{\beta}$ since, as we observed above, it is orthogonal to $Z(\mathfrak{h})$ and $ker(\chi) \subseteq Z(\mathfrak{h})$). Finally, it is easy to verify that $B(Z_F^0, Z_F^0) = -1$.

\smallskip

\noindent If $m=1$, then by $Z({\mathfrak{l}}) = \ker \chi$ and the orthogonality condition we immediately see that  $Z^0_F = \frac{Z^{\chi}}{\sqrt{-\| Z^{\chi}\|^2}}$

\medskip

\noindent We now calculate $\kappa$ defined by (\ref{DefKappa}).

\noindent To this aim, recall that $L$ is the isotropy subgroup of a non-zero vector $e \in \bC^m$ (with respect to the action of $H$ on $\bC^m$ defined through $\chi$).

\noindent In the case $m>1$, let us consider just the case when $\beta = \alpha_1$, $e=e_1$ and $Z^{\beta}$ is given by (\ref{Zbeta1}) (the calculation being similar in the case $\beta = \alpha_{m-1}$, $e = e_m$).

\noindent Since through the (Lie algebra) representation $t Z^{\chi}$ corresponds to $it Id_m$ and $t Z^{\beta}$ to $i t \ diag(m-1, -1, \dots, -1)$, by (\ref{explicitZ^0}) we have that $t Z^0_F$ goes in the Lie algebra representation to

$$\frac{1}{\sqrt{-\frac{1}{\| Z^{\chi}\|^2} - \frac{(m-1)^2}{\| Z^{\beta}\|^2}}} \left( \frac{it}{\| Z^{\chi}\|^2} Id + \frac{m-1}{\| Z^{\beta}\|^2} it \ diag(m-1, \dots ) \right) = $$

$$= diag \left( - it \sqrt{-\frac{1}{\| Z^{\chi}\|^2} - \frac{(m-1)^2}{\| Z^{\beta}\|^2}}, \dots \right)$$

\noindent so that $exp(tZ^0_F)$ goes to

$$diag \left( e^{-it \sqrt{-\frac{1}{\| Z^{\chi}\|^2} - \frac{(m-1)^2}{\| Z^{\beta}\|^2}} }, \dots \right).$$

\noindent So, in order for $exp(tZ^0)$ to fix $e_1$ we must have $-t \sqrt{-\frac{1}{\| Z^{\chi}\|^2} - \frac{(m-1)^2}{\| Z^{\beta}\|^2}} = 2 \pi k$ for some $k \in \bZ$, and the first positive value for which this holds true is $T_ 0 = \frac{2 \pi}{\sqrt{-\frac{1}{\| Z^{\chi}\|^2} - \frac{(m-1)^2}{\| Z^{\beta}\|^2}} }$, from which we finally deduce that

$$\kappa = \frac{2 \pi}{T_0} = \sqrt{-\frac{1}{\| Z^{\chi}\|^2} - \frac{(m-1)^2}{\| Z^{\beta}\|^2}} $$

\noindent that is

$$\kappa Z^0_F =   \frac{Z^{\chi}}{\| Z^{\chi}\|^2} + (m-1) \frac{Z^{\beta}}{\| Z^{\beta}\|^2}$$

\noindent (notice that this equality holds true both when $\beta = \alpha_1$ and $\beta = \alpha_{m-1}$).

\noindent Moreover, by the definition of $B(X,Y)$ in each of the classical groups recalled above and by (\ref{Zbeta1imm})-(\ref{Zbeta32}) we have $\| Z_{\beta} \|^2 = -2c m(m-1)$ where $c =  n, 2(n+1), 2(n-1), 2n-1$ for $\mathfrak{g} =  \mathfrak{su}_n,  \mathfrak{sp}_{n},  \mathfrak{so}_{2n},  \mathfrak{so}_{2n+1}$ respectively, so we finally get

\begin{equation}\label{explicitKZ^0}
\kappa Z^0_F =   \frac{Z^{\chi}}{\| Z^{\chi}\|^2} - \frac{Z^{\beta}}{2cm}
\end{equation}

\smallskip

\noindent In the case $m=1$, the same argument shows that $\kappa = \frac{1}{\sqrt{-\| Z^{\chi}\|^2}}$, so that we get

\begin{equation}\label{explicitKZ^0m=1}
\kappa Z^0_F =   -\frac{Z^{\chi}}{\| Z^{\chi}\|^2}
\end{equation}

\medskip

\noindent Now, we are going to rewrite this by $B$-duality, i.e. to calculate the dual form $\xi^0 = B^{-1}(\kappa Z^0_F)$. For the sake of brevity, from now on we will denote $Z \simeq \xi$ to mean $Z = B^{-1}(\xi)$.

\smallskip

\noindent First, the fact that $B(Z^{\chi}, Z) = 0$ for every $Z \in ker(\chi)$ means that $Z^{\chi} \simeq D \chi$, for some $D \in \bC$; then, by $\|Z^{\chi}\|^2 = B(Z^{\chi}, Z^{\chi}) = D \chi(Z^{\chi}) = D \ i$ we have

\begin{equation}\label{Zchifrattosimeq}
\frac{Z^{\chi}}{\| Z^{\chi}\|^2} \simeq \frac{\chi}{i}
\end{equation}

\noindent By (\ref{explicitKZ^0m=1}), this immediately yields $\kappa Z^0_F \simeq i \chi$ in the case $m=1$.

\medskip

\noindent In the case $m >1$, we need to calculate $B^{-1}(Z^{\beta})$.

\noindent Recall that, if we denote by $E_{ij}$ the square matrix having $1$ at position $i \ j$ and zero otherwise, then

$$\frac{1}{2c} (E_{ii} - E_{jj}) \simeq \varepsilon_i - \varepsilon_j \ \rm{for} \ \ \mathfrak{g} = \mathfrak{su}_{n} \ \ (c=n)$$

$$\frac{1}{2c} \left( \begin{array}{cc}
E_{ii} - E_{jj} & 0 \\
0 & E_{jj} - E_{ii}
\end{array} \right) \simeq \varepsilon_i - \varepsilon_j \ \rm{for} \ \ \mathfrak{g} = \mathfrak{sp}_{2n}, \mathfrak{so}_{2n} \ \ (c = 2(n+1), 2(n-1)  \ \rm{respectively})$$

$$\frac{1}{2c} \left( \begin{array}{ccc}
E_{ii} - E_{jj} & 0 & 0 \\
0 & E_{jj} - E_{ii} & 0 \\
0 & 0 & 0
\end{array} \right) \simeq \varepsilon_i - \varepsilon_j \ \rm{for} \ \ \mathfrak{g} = \mathfrak{so}_{2n+1}  \ \ (c= 2n-1)$$

\noindent Then, if $\alpha_1 = \varepsilon_k - \varepsilon_{k+1}$ and $\alpha_{m-1} = \varepsilon_{k+m-2} - \varepsilon_{k+m-1}$,  combining these identities with  (\ref{Zbeta1imm})-(\ref{Zbeta32}) above we get

\begin{equation}\label{ZbetaSimeq1}
Z^{\beta} \simeq 2c \ i ((m-1) \varepsilon_k - (\varepsilon_{k+1}  + \cdots + \varepsilon_{k+m-1}))
\end{equation}

\begin{equation}\label{ZbetaSimeq2}
( \ {\rm resp.} \ Z^{\beta} \simeq 2c \ i ((m-1) \varepsilon_{k+m-1} - (\varepsilon_{k} + \cdots + \varepsilon_{k+m-2}))  \ )
\end{equation}

\noindent for $\beta = \alpha_1$ (resp. $\beta = \alpha_{m-1}$).

\noindent Now, if we are not in one of the cases (1) or (2) of the statement of the lemma, then $\beta_s = \varepsilon_{k-1} - \varepsilon_{k}$, $\beta_{s+1}= \varepsilon_{k+m-1} - \varepsilon_{k+m}$ and the corresponding fundamental weights are

$$\pi_{s} = \varepsilon_1 + \cdots + \varepsilon_{k-1}, \ \ \ \pi_{s+1} = \varepsilon_1 + \cdots + \varepsilon_{k+m-1}$$

Moreover, the fundamental weight associated to the new black node $\beta = \alpha_1$ (resp. $\beta = \alpha_{m-1}$) is given by $\pi_0 =\varepsilon_1 + \cdots + \varepsilon_{k}$ (resp. $\pi_0 =\varepsilon_1 + \cdots + \varepsilon_{k+m-2}$) and then, by (\ref{ZbetaSimeq1}) and (\ref{ZbetaSimeq2}) we get

\begin{equation}\label{ZbetaSimeq1fw}
Z^{\beta} \simeq 2c \ i (m \pi_0 - (m-1) \pi_{s} - \pi_{s+1})
\end{equation}

\begin{equation}\label{ZbetaSimeq2fw}
( \ {\rm resp.} \ Z^{\beta} \simeq 2c \ i ((m-1) \pi_{s+1} - m \pi_0 + \pi_{s}) \ )
\end{equation}

\noindent If we are in case (1) of the statement of the lemma, we have $\pi_{s+1} = \frac{\epsilon_1 + \cdots + \epsilon_n}{2}$ and then

\begin{equation}\label{ZbetaSimeq1fwSO2n+1}
Z^{\beta} \simeq 2c \ i (m \pi_0 - (m-1) \pi_{s} -2 \pi_{s+1})
\end{equation}

\begin{equation}\label{ZbetaSimeq2fwSO2n+1}
( \ {\rm resp.} \ Z^{\beta} \simeq 2c \ i (2(m-1) \pi_{s+1} - m \pi_0 + \pi_{s}) \ )
\end{equation}

\noindent for $\beta = \alpha_1$ (resp. $\beta = \alpha_{m-1}$).

\noindent If we are in case (2) of the statement of the lemma, then $\pi_{s+1} = \frac{\epsilon_1 + \cdots + \epsilon_{n-1} - \epsilon_n}{2}$ and one sees that

\begin{equation}\label{ZbetaSimeq1fwSO2n}
Z^{\beta} \simeq 2c \ i (m \pi_0 - (m-1) \pi_{s} -2 \pi_{s+1})
\end{equation}

\begin{equation}\label{ZbetaSimeq2fwSO2n}
( \ {\rm resp.} \ Z^{\beta} \simeq 2c \ i ((m-2) \pi_{s+1} - m \pi_0 + \pi_{s}) \ )
\end{equation}

\noindent Then the lemma follows by substitution in  (\ref{explicitKZ^0}) and by $Z^0 = - i Z^0_F$ (recall also that both $Z^0$ and $Z^0_F$ are determined up to the sign).
\qed

\bigskip

Now we are ready to prove Theorems \ref{prop triple2} and \ref{prop triple4}. Let us recall that, assuming that the group $G$ is simply connected, it is known that the lattice of characters coincide with the lattice of weights (\cite{G-O-V}, \cite{Sn}) so that the  Lie algebra character $\chi$ is given by a linear combination {\it with integer coefficients} of the fundamental weights $\pi_1, \dots, \pi_p$ associated to the black nodes of the diagram of $G/H$:

$$\chi = \sum_{j=1}^p k_j \pi_j, \ \ k_j \in \bZ$$

\noindent while the Koszul form $\sigma = B^{-1}(Z^{Kos}) = \sum_{\alpha \in R^+_{\mathfrak{m}}}\alpha$ (where
$R^+_{\mathfrak{m}} =  R^+ \setminus R^+_{\mathfrak{h}}$ denotes the set of  complementary  to $R^+_{\mathfrak{h}}$ positive  roots). In what follows, when necessary to avoid ambiguity we will write $Z^{Kos}_{G/H}$ to denote the Koszul vector of the flag manifold $G/H$.

\medskip

\noindent {\bf Proof of Theorem \ref{prop triple2}: \ }

Let $\pi_1, \dots, \pi_p$ be the fundamental weights of the black roots $\beta_1, \dots, \beta_p$ of the Dynkin diagram of $G/H$, and $n_1, \dots, n_p$ be the corresponding Koszul numbers, so that $B^{-1}(Z^{Kos}_{G/H}) = \sum_{j=1}^p n_j \pi_j$.

\noindent Let $m > 1$, let $A_{m-1} = \{ \alpha_1, \dots, \alpha_{m-1} \}$ be the white string given by the data which define the admissible vector bundle and let $\pi_0$ be the fundamental weight of the new black node $\beta$ (with $\beta = \alpha_1$ or $\beta = \alpha_{m-1}$), so that we have $B^{-1}(Z^{Kos}_{G/K}) = n_0' \pi_0 + \sum_{j=1}^p n_j' \pi_j$.

\noindent By using the description of the Koszul numbers given in subsection  2.6, one easily verifies that if $\beta = \alpha_1$ (resp. $\beta = \alpha_{m-1}$), then

\begin{equation}\label{nsn'sSU}
n_ 0' = m, \ \  n'_s = n_s - (m -1), \ \ n_{s+1}' =  n_{s+1} - 1
\end{equation}

\begin{equation}\label{nsn'sSU2}
{\rm (resp. } \ \ n_ 0' = m, \ \  n'_s = n_s - 1, \ \ n_{s+1}' =  n_{s+1} - (m-1) \ \ )
\end{equation}

\noindent with the exception of the same two cases seen in Lemma \ref{LemmaKZ^0}, that is

\begin{equation}\label{paintedG/H1}
\mathfrak{g} = \mathfrak{so}_{2n+1}, \ \ \ \cdots -  \underset{\beta_s}\bullet - \underset{\alpha_1}\circ - \cdots - \underset{\alpha_{m-1}}\circ \Rightarrow \underset{\beta_{s+1}}\bullet
\end{equation}

\noindent where, if $\beta = \alpha_1$ (resp. $\beta = \alpha_{m-1}$), we have

\begin{equation}\label{nsn'sSOdispari}
n_ 0' = m, \ \  n'_s = n_s - (m-1), \ \ n_{s+1}' =  n_{s+1} - 2
\end{equation}

\begin{equation}\label{nsn'sSOdispari2}
{\rm (resp.} \ \ n_ 0' = m, \ \  n'_s = n_s - 1, \ \ n_{s+1}' =  n_{s+1} - 2(m-1) \ )
\end{equation}

\noindent and

\begin{align}\label{paintedG/H2}
\mathfrak{g} = \mathfrak{so}_{2n}, \ \ \  \cdots - \underset{{\beta_s}}{\bullet} - \underset{{\alpha_1}}{\circ} - \dotsb - \underset{{}}{\overset{\overset{\textstyle\bullet_{{\beta{s+1}}}}{\textstyle\vert}}{\circ}} \,-\, \underset{{\alpha_{m-1}}}{\circ} \\
\end{align}

\noindent where

\begin{equation}\label{nsn'sSOpari}
n_ 0' = m, \ \  n'_s = n_s - (m-1), \ \ n_{s+1}' =  n_{s+1} - 2
\end{equation}

\begin{equation}\label{nsn'sSOpari2}
{\rm (resp.} \ \ n_ 0' = m, \ \  n'_s = n_s - 1, \ \ n_{s+1}' =  n_{s+1} - (m-2) \ )
\end{equation}

\noindent for $\beta = \alpha_1$ (resp. $\beta = \alpha_{m-1}$).

\noindent We point out that formulas (\ref{nsn'sSU}) and (\ref{nsn'sSU2}) hold true also for $G = Sp(n)$ in the case when the painted Dynkin diagram of $G/H$ is

\begin{equation}\label{paintedG/H}
\cdots -  \underset{\beta_s}\bullet - \underset{\alpha_1}\circ - \cdots - \underset{\alpha_{m-1}}\circ \Leftarrow \underset{\beta_{s+1}}\bullet
\end{equation}

\bigskip

\noindent Now, if $\lambda=0$, the algebraic condition for the existence of the Einstein metric is

\begin{equation}\label{condizlambda0prima}
Z^{Kos} = \kappa m Z^0
\end{equation}

By using Lemma \ref{LemmaKZ^0} and (\ref{nsn'sSU}),(\ref{nsn'sSU2}), (\ref{nsn'sSOdispari}), (\ref{nsn'sSOdispari2}), (\ref{nsn'sSOpari}), (\ref{nsn'sSOpari2}) if $\beta = \alpha_1$ (resp. $\beta = \alpha_{m-1}$), one rewrites condition (\ref{condizlambda0prima}) as

\begin{equation}\label{ZkosZ01}
m \pi_0 + \sum_{j \neq s, s+1} n_j \pi_j + [n_s-(m-1)]\pi_s + [n_{s+1} - 1] \pi_{s+1} = m \chi + m \pi_0 - (m-1) \pi_s - \pi_{s+1}
\end{equation}

\noindent (resp.

\begin{equation}\label{ZkosZ01resp}
m \pi_0 + \sum_{j \neq s, s+1} n_j \pi_j + (n_s-1)\pi_s + [n_{s+1} - (m-1)] \pi_{s+1} = m \pi_0 - m \chi - (m-1) \pi_{s+1} - \pi_{s} \ )
\end{equation}

\noindent or

\begin{equation}\label{ZkosZ01SO2n+1}
m \pi_0 + \sum_{j \neq s, s+1} n_j \pi_j + [n_s-(m-1)]\pi_s + [n_{s+1} - 2] \pi_{s+1} = m \chi + m \pi_0 - (m-1) \pi_s - 2\pi_{s+1}
\end{equation}

\noindent (resp.

\begin{equation}\label{ZkosZ01SO2n+1resp}
m \pi_0 + \sum_{j \neq s, s+1} n_j \pi_j + (n_s-1)\pi_s + [n_{s+1} - 2(m-1)] \pi_{s+1} = m \pi_0 - m \chi - 2(m-1) \pi_{s+1} - \pi_{s} \ )
\end{equation}

\smallskip

\begin{equation}\label{ZkosZ01SO2n}
m \pi_0 + \sum_{j \neq s, s+1} n_j \pi_j + [n_s-(m-1)]\pi_s + [n_{s+1} - 2] \pi_{s+1} = m \chi + m \pi_0 - (m-1) \pi_s - 2\pi_{s+1}
\end{equation}

\noindent (resp.

\begin{equation}\label{ZkosZ01SO2nresp}
m \pi_0 + \sum_{j \neq s, s+1} n_j \pi_j + (n_s-1)\pi_s + [n_{s+1} - (m-2)] \pi_{s+1} = m \pi_0 - m \chi - (m-2) \pi_{s+1} - \pi_{s} \ )
\end{equation}

in the exceptional cases of $\mathfrak{g} = \mathfrak{so}_{2n+1}$ and $\mathfrak{g} = \mathfrak{so}_{2n}$.

In all the cases (\ref{ZkosZ01}), (\ref{ZkosZ01resp}), (\ref{ZkosZ01SO2n+1}), (\ref{ZkosZ01SO2n+1resp}), (\ref{ZkosZ01SO2n}), (\ref{ZkosZ01SO2nresp}), one immediately sees that after simplifications one gets

$$\sum_{j = 1}^p n_j \pi_j = m \chi$$

\noindent from which (i) of Theorem \ref{prop triple2} follows.

\medskip

In order to prove (ii), recall that in the case $\lambda \neq 0$, the necessary and sufficient conditions to have a standard Kahler-Einstein metric in a neighbourhood of the singular section are (see Theorem \ref{maintheoronKE})

\begin{enumerate}
\item[(1)]  the vector $Z_0 = \frac{1}{\lambda} (Z^{Kos} - km Z^0)$ \noindent satisfies $\beta_j(Z_0) > 0$, $j=1, \dots, p$ and $\beta(Z_0) = 0$
\item[(2)] for at least small values of $s>0$, the segment $Z_0 + s Z^0$ satisfies $\beta_j(Z_0 + s Z^0) > 0$, $j=1, \dots, p$ and  $\beta(Z_0 + s Z^0) > 0$
\end{enumerate}

\noindent In fact, if (1) is satisfied, then $\beta(Z_0 + s Z^0) >0$ reduces to $\beta(Z^0) >0$ and $\beta_j(Z_0 + s Z^0) > 0$ is always satisfied for small values of $s$.

\noindent So we need just to check $\beta(Z^0) > 0$, $\beta_j(Z_0) > 0$, $\beta(Z_0) = 0$.

\smallskip

\noindent The first condition is easily verified by (\ref{KZ^0INLEMMA22}), (\ref{KZ^0INLEMMA24lemma}), (\ref{KZ^0INLEMMA23lemma}) in Lemma  \ref{LemmaKZ^0}.

As for the other two conditions, by (\ref{ZkosZ01}), (\ref{ZkosZ01resp}), (\ref{ZkosZ01SO2n+1}), (\ref{ZkosZ01SO2n+1resp}), (\ref{ZkosZ01SO2n}), (\ref{ZkosZ01SO2nresp}), one sees immediately that in any case

$$B^{-1} (Z_0) = \frac{1}{\lambda}(\sum_{j=1}^p n_j \pi_j - m \chi) $$

\begin{equation}\label{Z0allcases}
({\rm resp.} \ \ B^{-1} (Z_0)  = \frac{1}{\lambda}(\sum_{j=1}^p n_j \pi_j + m \chi) \ )
\end{equation}

\noindent for $\beta = \a_1$ (resp. $\beta = \a_{m-1}$), so $\beta(Z_0) = 0$ is clear, while $\beta_j(Z_0) > 0$  writes

\begin{equation}\label{conditionlambdadiverso042}
\left. \begin{array}{c}
\frac{1}{\lambda} \frac{\| \beta_j \|^2}{2}  [n_j - m k_j ] > 0, \ \ j = 1, \dots, p
\end{array} \right.
\end{equation}

\begin{equation}\label{conditionlambdadiverso042}
\left. \begin{array}{c}
{\rm (resp. \ \ } \frac{1}{\lambda} \frac{\| \beta_j \|^2}{2}  [n_j + m k_j ] > 0, \ \ j = 1, \dots, p
\end{array} \right.
\end{equation}

\noindent from which the assertions in (ii) of Theorem \ref{prop triple2} immediately follow. \qed

\bigskip

\noindent {\bf Proof of Theorem \ref{prop triple4}: \ } The proof is the same as for Theorem \ref{prop triple2}, but now  since $m=1$ we have $G/K = G/H$ and the Koszul forms $\sigma$ of $G/K$ and $G/H$ coincide: then, for $\lambda = 0$ we use (\ref{KZ^0INLEMMA1}) to rewrite condition (\ref{condizlambda0prima}) for the existence of a Ricci-flat metric as

\begin{equation}\label{condizlambda0}
\chi = \sigma = \sum_{j=1}^p n_j \pi_j
\end{equation}

\noindent while, in the case $\lambda \neq 0$, since there is no new black root $\beta$ and the only condition to have a standard K\"ahler-Einstein metric in a neighbourhood of the singular section is $\beta_j (Z_0)>0$, the assertions of the theorem immediately follow from $B^{-1} (Z_0) = \frac{1}{\lambda}(\sum_{j=1}^p n_j \pi_j - \chi)$. \qed

\begin{remar}\label{remunique}
\rm
The above proofs show that, once the conditions for the existence of a K\"ahler-Einstein standard \C1 structure on an admissible vector bundle $M_{\vp} =  G \times_H V_{\vp}$ having $S_0$ as only singular orbit are satisfied, then this structure is unique provided the Einstein constant $\lambda \neq 0$, while for $\lambda =0$ the metrics are parametrized by the vectors $Z_0 \in C(J^S)$. Indeed, in the case $\lambda =0$ condition (\ref{condizlambda0prima}) for the existence of the metric does not depend on the choice of  the initial vector $Z_0 \in C(J^S)$ of the segment in $C(J^F)$, while for $\lambda \neq 0$ the vector $Z_0$ is completely determined by $Z_0 = \frac{1}{\lambda}(Z^{Kos} - km Z^0)$. This is consistent with the results in \cite{D-W}.
\end{remar}

\begin{remar}\label{RemarkKpos}
\rm
By the same remarks made to prove (ii) of Theorem \ref{prop triple2}, one sees that, if a metric with $\lambda \neq 0$ exists, then the segment $Z_0 + x Z^0$ can be extended to a whole ray in the T-Weyl chamber $C(J^F)$ of $F = G/K$ if and only if $\beta_j(Z^0) > 0$ for every $j=1, \dots, p$.

\noindent If we are not in one of the exceptional cases of the statement of Lemma \ref{LemmaKZ^0}, by (\ref{KZ^0INLEMMA22}) this condition both for $\lambda > 0$ and $\lambda < 0$ reads

$$k_j > 0, \ \   k_s > \frac{m-1}{m}, \ \ k_{s+1} > \frac{1}{m} $$
\begin{equation}\label{condizBordo1}
( \ {\rm resp.} \ k_j < 0, \ \   k_s < -\frac{1}{m}, \ \ k_{s+1} < -\frac{m-1}{m} \ )
\end{equation}

\noindent if $\beta = \a_1$ (resp. $\beta = \a_{m-1}$).

\noindent By Theorem \ref{prop triple2}, these conditions are always compatible with those for the existence of the K\"ahler-Einstein metric in the case $\lambda < 0$, i.e. $k_j > \frac{n_j}{m}$ (resp. $k_j < - \frac{n_j}{m}$), while this is not true in the case $\lambda > 0$, where we must find {\it integer} $k_j$'s, $j=1, \dots, p$, such that $0 < k_j < \frac{n_j}{m}$ for $j \neq s, s+1$ and $\frac{m-1}{m} < k_s < \frac{n_s}{m}$, $\frac{1}{m} < k_{s+1} < \frac{n_{s+1}}{m}$ (resp. $- \frac{n_j}{m} < k_j < 0$ for $j \neq s, s+1$ and $ -\frac{n_s}{m} < k_s < - \frac{1}{m}$, $ -\frac{n_{s+1}}{m} < k_{s+1} < - \frac{m-1}{m}$).

\noindent In the case $\lambda = 0$, the algebraic condition (\ref{condizlambda0prima}) implies that $Z^0 = \frac{1}{\kappa m} Z^{Kos} \in C(J^F)$, and then for any choice of the starting point $Z_0 \in C(J^F) \cap \{ \beta=0 \}$, the segment $Z_0 + x Z^0$ extends to a ray in $C(J^F)$.

\noindent This is in accordance with the last assertion in Theorem \ref{maintheoronKE} (see also the end of the proof of Theorem 36 in \cite{AZ}).

\noindent We leave the details about the exceptional cases of the statement of Lemma \ref{LemmaKZ^0} to the reader.

\end{remar}

\begin{remar}
\rm Notice that $k_j \in \bZ$ assures that the vector $Z^0$ fulfills the condition that $\{ \exp(t Z^0) \ | \ t \in \bR\}$ is compact (more precisely, a circle $S^1$)
\end{remar}

\begin{remar}
\rm Under the conditions given by the proposition, there always exists at least a non-complete Kahler-Einstein metric, given by a segment $Z_0 + s Z^0$ staying in the interior of the chamber $C$ and parametrized with $s = f(t)$, being $f(t)$ the solution to the ODE (\ref{EinstEq}).

\noindent As we have proved in the first part of this paper, if the conditions of Proposition \ref{prop triple2} are satisfied, there exists a complete Kahler-Einstein metric for $\lambda \leq 0$, while for $\lambda > 0$ the metric is never complete.

\end{remar}

\appendix
\section{ Appendix. Basic  facts  on flag manifolds}

\noindent  Let $F = G/K = \mathrm{Ad}_G Z$, where $Z \in
\mathfrak{g}$, be a  flag manifold, i.e. an adjoint orbit of a
 compact semisimple Lie group $G$
  with   the  $B$-orthogonal \\
   (where $B$ is  the Killing form) reductive   decomposition
  $$ \mathfrak{g} = \gk + \gm  = C_{\gg}(Z) + \gm.$$
   We can decompose  $\gk$  as
     $$\gk   = Z(\gk) \oplus \gk'$$
   where  $\gk'$ is  the semisimple part   and  $Z(\gk)$ is  the center.
    We fix a Cartan subalgebra  $\gc$ of  $\gk$ (hence  also of
    $\gg$) and   denote  by
      $R$  the   root  system
     of   the   complex Lie  algebra $\gg^{\mathbb{C}}$ w.r.t. the Cartan  subalgebra
      $\gc^{\mathbb{C}}$.  We set
      $$R_{\mathfrak{k}} :=\{ \a \in R, \, \a(Z(\gk)) = 0 \},\,R_{\mathfrak{m}} := R \setminus R_{\mathfrak{k}}. $$

 Then
  $$ \gk = \gc  + \gg (R_{\mathfrak{k}})^{\tau},\,   \gm = \gg(R_{\mathfrak{m}})^{\tau}, $$
where for a  subset $P \subset R$, we set
$$\gg(P) = \sum_{\alpha \in
P} \gg_{\alpha}$$

\noindent being $\gg_{\a}$ the  root  space with root  $\a$ and
$V^{\tau}$ means the fix point set in $V\subset \gg^{\bC}$  of the
complex conjugation $\tau$. Recall that the Killing   form induces
an Euclidean metric in the real vector space $i \gc $ and roots are
identified  with real linear forms on $i \gc$. We  set $\gt:= i
Z(\gk) \subset i \gc$ and denote  by
 $$\rho: R \to R|_{\gt},\,\,  \a \mapsto \bar \a := \a|_{\gt}$$
  the  restriction map.

\begin{definition}
The set $R_T =\rho(R_{\mathfrak{m}})= R_{\mathfrak{m}}|{\gt}$ of linear forms
on $\gt$ which are  restriction of  roots from $R_{\mathfrak{m}}$ is called the system of
{\bf $T$-roots} and connected components $C$ of the set $\gt
\setminus \{ \ker \bar \a ,\, \bar \a \in R_T \}$
  are called  {\bf $T$-Weyl chambers}.
\end{definition}

Sets of  $T$-roots $\xi$  bijectively correspond to irreducible   $\gk$-submodules
 $ \gm(\xi):=  \gg(\rho^{-1}(\xi))$ of the   complexified  isotropy module
 $\gm^{\mathbb{C}}$ of  the  flag manifold $F =G/K$.

So a decomposition of the $\gk$-modules  $\gm^{\mathbb{C}}$ and $\gm$
into irreducible
  submodules   can be  written
 as
$$  \gm^{\mathbb{C}} = \sum_{\xi \in R_T} \gm(\xi),\,\,
  \gm = \sum_{\xi \in R^+_T}[ \gm(\xi) + \gm(-\xi) ]^{\tau}$$
 where $R^+_T :=\rho(R_{\mathfrak{m}}^+)$ is   the  system of  positive $T$-roots associated    with  a
system of positive  roots $R^+$, see  \cite{A-P}, \cite{A}.

We  fix a system of  simple  roots $\Pi_W$ of  $R_{\mathfrak{k}}$
and denote  by  $\Pi = \Pi_W \cup \Pi_B$ its extension  to a system
of simple roots of $R$.  Let $R^+ = R^+(\Pi) $  be the associated
system of positive  roots and $R^+_{\mathfrak{m}}:= R^+ \cap
R_{\mathfrak{m}}$. The  set  $R^+_T:= \rho(R^+_{\mathfrak{m}})$ is
called positive $T$-root  set.

 We need  the following

\begin{theorem}\cite{A-P}
There exists a one-to-one correspondence  between
extensions $\Pi = \Pi_W \cup \Pi_B$  of  the system $\Pi_W$ of
simple system of  $R_{\mathfrak{k}}$,  $T$-Weyl chambers $C\subset
\gt$ and invariant complex structures (ICS) $J$  on $F = G/K$. If
$\Pi_B = \{ \beta_1, \dots ,\beta_k \}$, then  the corresponding
$T$-Weyl chamber  is defined  by  $   C = \{   \bar \beta_1>0,
\dots ,\bar \beta_k >0 \}$ where $\bar \beta = \rho(\beta)$  and
the  complex  structure is  defined  by $\pm i$-eigenspace decomposition
 \be \label{InvCompStrin"m"} \gm^{\bC} = \gm^+ + \gm^- =
\gg(R^+_{\mathfrak{m}}) + \gg(-R^+_{\mathfrak{m}}) \ee
 of the
complexified  tangent  space  $\gm^{\bC} = T_{eK}(G/K)$.
 \end{theorem}

The extension $\Pi = \Pi_W \cup \Pi_B$  can be graphically described
by  a painted Dynkin diagram, i.e.  the  Dynkin  diagram  which
represents the  system $\Pi$  with the nodes representing  $\Pi_B$
painted in  black. Such a diagram, which we   sometimes  identify
with  the pair $(\Pi_W, \Pi_B)$,  allows  to reconstruct   the flag
manifold $F= G/K$ with invariant  complex  structure $J^F$ as
follows: the  semisimple  part $\gk'$ of the  (connected) stability
subalgebra $\gk$ is  defined as the  regular semisimple  subalgebra
associated   with the closed  subsystem $R_{\mathfrak{k}}= R \cap
\mathrm{span}(\Pi_W)$ and the vectors $ih_{j}$ defined  by condition
$$ \beta_k(h_j) = \delta_{kj}, \, \alpha_i(h_j)=0,\,  \beta_j \in \Pi_B, \alpha_i \in \Pi_W$$
 form a basis of
the center $Z(\mathfrak{k})$. The complex structure is defined  by
(\ref{InvCompStrin"m"}).

\smallskip

\noindent Now, an element $Z \in \mathfrak{t}$  is called   to  be {\bf
$K$-regular} if  its   centralizer  $C_G(Z) = K$ or, equivalently,
any $T$-root has a non-zero value on $Z$. Then we have the following

\begin{proposition}\label{propcorrespmetrics}
(\cite{B-H}, \cite{A-P}) There exists
a natural one-to-one correspondence between elements $Z \in
\mathfrak{t}$ and closed invariant 2-forms  $\omega_Z$ on $G/K$,
given by
$$Z \leftrightarrow \omega_Z|_o = i \, d(B \circ Z),$$
\noindent where $d$
is  the exterior differential in  the Lie  algebra $\mathfrak{g}$
defined  by $d\alpha(X,Y) = - 1/2\alpha([X,Y])$  and  $o = eK \in G/K$.\\
Moreover,  regular elements $Z \in C$  from a $T$-Weyl chamber  $C$
correspond to the  K\"ahler  forms  $\omega_Z$ with respect to the
complex  structure   $J(C)$ associated to $C$, that is they define an invariant  K\"ahler
structure  $(\omega_Z, J(C))$.  The  2-form $ \frac{1}{2\pi}\omega_Z$
is integral if  the  1-form  $B\circ Z$  has integer  coordinates
 with respect  to  the  fundamental weights  $ \pi_i$   associated  with
  the  system  of black simple roots $ \beta_i \in \Pi_B$.
\end{proposition}

\noindent Recall that if $\Pi_W = \{ \a_1, \dots, \a_m\}$ (resp. $\Pi_B = \{ \beta_1, \dots, \beta_k\}$) is the set of
 white
(resp. black)  simple roots, then the fundamental weight $ \pi_i$
associated with $\beta_i$, $i = 1, \dots, k$, is the linear form
defined by

\begin{equation}\label{deffundwe}
\frac{2\langle  \pi_i, \beta_j \rangle}{\| \beta_j \|^2} =
\delta_{ij}, \ \ \langle  \pi_i, \alpha_j \rangle = 0.
\end{equation}
where $< .,.>$ is the scalar product  in $i\gc^* = \mathrm{span}(R)$
induced by the Killing  form.
\noindent The $B$-dual  to $ \pi_i$ vectors $h_i$ form
  a basis of $\gt$.

\noindent Let $E_{\alpha} \in \gg_{\alpha}, \, \alpha \in R$, be the {\it Chevalley basis} of $\mathfrak{g}(R)$
such that $B(E_{\alpha}, E_{- \alpha}) = \frac{2}{<\alpha, \alpha>}$
  We denote by   $\omega_{\alpha} =
B\circ E_{\alpha}$ the dual basis of 1-forms. Then for $Z \in \gt$

\be \label{omega_Z}
 \omega_Z = -i\sum_{\alpha \in R_{\mathfrak{m}}^+}
\frac{2 \alpha(Z)}{<\alpha, \alpha>} \omega_{\alpha} \wedge
\omega_{-\alpha}\ee

\noindent Indeed,

$$
\begin{array}{ll} i \, d(B\circ Z)(E_{\alpha}, E_{- \alpha})  =&
-\frac{i}{2} B(Z, [E_{\alpha},
E_{-\alpha}])\\
=&-\frac{i}{2}B([Z,E_{\alpha}],E_{-\alpha})\\
= &-  \frac{i}{2}\alpha(Z)B(E_{\alpha}, E_{-\alpha})\\
=&- \frac{i\alpha(Z)}{<\alpha, \alpha>}\\
 = &-2i\frac{\alpha(Z)}{<\alpha, \alpha>} \omega_{\alpha} \wedge
\omega_{-\alpha}(E_{\alpha}, E_{-\alpha}).
\end{array}
$$

\begin{definition}\label{DefKosVect}
The 1-form
$$\sigma = \sum_{\beta \in
R_{\mathfrak{m}}^+} \beta \in \gt^* \subset i\gc^*$$
 is called  the
Koszul form and  the dual vector $Z^{Kos}:= B^{-1} \circ \sigma$ is
called the Koszul vector.
\end{definition}

\begin{proposition}\label{propKoszul} \cite{A-P}
The {\bf Koszul vector}
 $Z^{Kos}$
 defines  the invariant  K\"ahler-Einstein
 structure  $(\omega_{Z^{Kos}}, J(C))$ on $F=G/K$, where
$J(C)$ is the invariant complex structure  associated  with the
$T$-Weyl chamber $C$ which is defined by $\Pi_B$.
\end{proposition}


\begin{thebibliography}{ABBR}

\bibitem{A} Alekseevsky D.: {\em Flag manifolds}, 11. Yugoslav Geometrical seminar, Divcibare, 10-17 October, 3-35 (1993)

\bibitem{AZ} Alekseevsky D., Zuddas F.: {\em Cohomogeneity one Kahler and Kahler-Einstein manifolds with one singular orbit I}, Ann. global Anal. Geom. 52:1, 99-128 (2017)
\bibitem{A-Ch} Alekseevsky D., Chrysikos J.:  {\em
Spin  structures  on  compact  homogeneous  pseudo-Riemannian manifolds}, Transf.  Groups (2018), pp. 1-31.




\bibitem{A-S} Alekseevsky D., Spiro A.: {\em Invariant CR structures on compact homogeneous manifolds},
 Hokk. Math. J, v. 32, no.2, 209-276 (2003)
\bibitem{A-P} Alekseevsky D. V., Perelomov, A. M.: {\em Invariant Kaehler-Einstein metrics on compact homogeneous
spaces}, Funct. Anal. Appl., 20 (3), 171-182 (1986)
 \bibitem{A-C-H=K} Alekseevsky D., Cortes V., Hasegawa K., Kamishima Y., {\em Homogeneous  locally  conformally K\"ahler  and Sasaki manifolds}, Int. J. Math.,  26, n5, ( 2015).
\bibitem{arv} Arvanitoyeorgos A.: {\em Geometry of flag manifolds}, International Journal of Geometric Methods in Modern Physics Vol.3, Nos. 5, 6, 957-974 (2006)

\bibitem{AZ-B} Achmed-Zade I., Bykov D.  , 
{\em Ricci-flat metrics on vector bundles over flag manifolds}, arXiv:1905.00412 (01/05/2019)

\bibitem{A-B} Azad H., Biswas I., {\em Quasi-potentials and Kahler Einstein metrics on  flag manifolds II}, Journal of Algebra 269 no. 2, (2003) 480–491.



\bibitem{B} Besse A.: {\em Einstein manifolds}, Ergeb. Math. Grenzgeb. (3) 10, Springer, Berlin, 1987.
\bibitem{Ber} Berard- Bergery L.: {\em Sur des nouvelles varietes Riemanniennes d'Einstein}, Publ.
de Inst. E. Cartan, No. 6, 1-60 (1982)
\bibitem{B-H} Borel A., Hirzebruch F.:  {\em Characteristic  classes
and  homogerneous  spaces}, Amer. J . Math. 80, 458-538 (1958)
\bibitem{B-G} Boyer C.P.,  Galicki K., {\em  Sasakian geometry}, Oxford Mathematical Monographs, Oxford University Press, Oxford, 2008.

\bibitem{C} Van Coevering C.,  {\em Calabi-Yau metrics on canonical bundles of flag varieties}, arXiv:1807.07256v1 (19/07/2018)

\bibitem{D-W} Dancer A., Wang M. Y.: {\em K\"ahler Einstein  metrics of
cohomogeneity one  and bundle construction for Einstein Hermitian
metrics}, Math. Ann. 312, 503-526 (1998)
\bibitem{E-W} Eschenburg J.-H., Wang M. Y.:
{\em The initial value problem for cohomogeneity one Einstein
metrics },  J. Geom. Anal.,10, No.1, 109-137 (2000)
\bibitem{G-O-V} Gorbatsevich V.V.,  Onishchik A.L., Vinberg E.B.: {\em Structure of  Lie  groups and Lie  algebras}, Encycl. Math. Sci.,
Lie groups  and Lie  algebras, III, Springer Verlag.
\bibitem{H-S} Huckleberry A., Snow D.: {\em Almost homogeneous K\"ahler manifolds with hypersurface orbits}, Osaka J.math. 19, 763-786 (1982)
\bibitem{K-S1}  Koiso N., Sakane Y.:
{\em Non-homogeneous K\"ahler-Einstein metrics on compact complex
manifolds}, Curvature and topology of Riemannian manifolds, Proc.
17th Int. Taniguchi Symp., Katata/Jap. 1985, Lect. Notes Math. 1201,
165-179 (1986)
\bibitem{K-S2} Koiso N., Sakane Y.: {\em Non homogeneous K\"ahler Einstein metrics on compact complex manifolds II}, Osaka J. Math. 25, 933-959 (1988)
\bibitem{P-P} Page D.N., Pope C.N.: {\em Inhomogeneous Einstein metrics on complex line bundles}, Classical Quantum Gravity 4 no. 2, 213-225 (1987)
\bibitem{P-S} Podest\`a F., Spiro A.: {\em Kaehler manifolds with large isometry group}, Osaka J. Math. Volume 36, Number 4, 805-833 (1999)
\bibitem{S} Sakane Y.: {\em Examples of compact Einstein-K\"ahler manifolds with positive Ricci tensor}, Osaka J. Math. 23, 585-616 (1986)
\bibitem{Sn} Snow D. M.: {\em Homogeneous vector bundles}, Group actions and invariant theory (Montreal,
PQ, 1988), CMS Conf. Proc., 10, 193-205, Amer. Math. Soc., Providence, RI (1989)
 \bibitem{Sp}  Sparks J.: {\em Sasaki-Einstein Manifolds}, Surveys Diff.Geom. 16, 265-324   (2011)
\bibitem{V} Verdiani L.: Invariant  metrics on cohomogeneity one manifolds,
Geometriae Dedicata ,  77 (1), 77-110 (1999)

\end{thebibliography}
\end{document}